\documentclass[11pt,reqno]{amsart}
\usepackage[T1]{fontenc}\usepackage{lmodern}
\usepackage[utf8]{inputenc}
\usepackage{amsmath,amssymb,mathrsfs}
\usepackage[left=0.95in,right=0.95in,top=0.95in,bottom=0.75in]{geometry}
\usepackage{microtype}
\usepackage[colorlinks=true,linkcolor=blue,citecolor=blue,urlcolor=blue]{hyperref}
\usepackage[
  backend=biber,
  style=trad-plain,
  url=false,
  doi=false,
  isbn=false
]{biblatex}
\DeclareFieldFormat[unpublished]{title}{#1}
\DeclareFieldFormat{eprint:arxiv}{arXiv\addcolon#1}

\DeclareBibliographyDriver{unpublished}{%
  \usebibmacro{bibindex}%
  \usebibmacro{begentry}%
  \printnames{author}%
  \newunit\newblock
  \usebibmacro{title}%
  \newunit\newblock
  \printtext{Preprint}%
  \setunit{\addcomma\space}%
  \usebibmacro{eprint}%
  \usebibmacro{finentry}%
}

\newcommand{\nocontentsline}[3]{}
\newcommand\stoptoc{%
   \let\origcontentsline\addcontentsline
   \let\addcontentsline\nocontentsline
}
\newcommand\resumetoc{%
   \let\addcontentsline\origcontentsline
}

\numberwithin{equation}{section}
\newtheorem{theorem}{Theorem}[section]
\newtheorem{proposition}[theorem]{Proposition}
\newtheorem{lemma}[theorem]{Lemma}
\newtheorem{corollary}[theorem]{Corollary}
\theoremstyle{remark}
\newtheorem{remark}[theorem]{Remark}
\newcommand{\R}{\mathbb R}
\newcommand{\C}{\mathbb C}
\newcommand{\Z}{\mathbb Z}
\newcommand{\RP}{\mathbb{RP}}
\newcommand{\CP}{\mathbb{CP}}
\newcommand{\Lt}{\widetilde L}
\newcommand{\Qt}{\widetilde Q}
\newcommand{\SSb}{\mathbb S}
\newcommand{\St}{\widetilde\Sigma}
\newcommand{\Fh}{\widehat F}
\newcommand{\Kc}{\mathcal K}
\newcommand{\HH}{\mathcal H}

\DeclareMathOperator{\Ind}{Ind}
\DeclareMathOperator{\Null}{Null}
\DeclareMathOperator{\tr}{tr}
\DeclareMathOperator{\spec}{spec}
\DeclareMathOperator{\Area}{Area}
\DeclareMathOperator{\Sym}{Sym}
\DeclareMathOperator{\Span}{span}
\DeclareMathOperator{\diver}{div}
\DeclareMathOperator{\ind}{ind}
\title[Morse index characterization of the Veronese plane]{%
  \mbox{Morse index characterization of the Veronese projective plane}}
\makeatletter
\renewcommand{\author@andify}{%
  \nxandlist{\unskip, }{\unskip, }{\unskip, }}
\makeatother
\author{Alessandro Carlotto}
\address{Universit\`a di Trento, Dipartimento di Matematica, via Sommarive 14, 38123 Trento, Italy}
\email{alessandro.carlotto@unitn.it}
\author{Yangyang Li}
\address{Department of Mathematics, University of Notre Dame, 218 Hayes-Healy Center, Notre Dame, IN 46556, United States of America}
\email{yangyangli@nd.edu}
\author{Zhihan Wang}
\address{Department of Mathematics, Cornell University, 310 Malott Hall, Ithaca, NY 14853, United States of America}
\email{zw782@cornell.edu}
\date{}
\AtBeginDocument{%
  \setlength{\abovedisplayskip}{5pt plus 2pt minus 2pt}%
  \setlength{\belowdisplayskip}{5pt plus 2pt minus 2pt}%
}
\begin{document}
\begin{abstract}
Following Simons' seminal work in 1968, every closed, connected, embedded minimal surface in the four-sphere that is not totally geodesic has Morse index at least five. In this paper, we prove that equality holds if and only if the surface is congruent to the Veronese projective plane. Relying on recent advances by Bernstein and Ketover, this characterization settles Kusner's 1996 conjecture on the minimizers of the Willmore functional in dimension four for surfaces whose normal bundle has nonzero Euler number. This is always the case, in particular, for embeddings of nonorientable surfaces whose orientation double cover has even genus.
\end{abstract}

\dedicatory{Dedicated to Professor Richard Schoen,
with gratitude and admiration,
on the occasion of his 76th birthday}

\renewcommand{\shortauthors}{Alessandro Carlotto, Yangyang Li, Zhihan Wang}

\maketitle
\setcounter{tocdepth}{1}
\tableofcontents
\section{Introduction}
Minimal surfaces are critical points of the area functional, and their Morse index  is the maximal dimension of a space of smooth normal variations on which the second variation is negative definite. The index thus quantifies instability. It also connects variational existence theorems to geometric classification: an upper bound on the index of a surface produced by a variational construction can become a rigidity statement when the surfaces of low index are understood.

A fundamental index estimate due to Simons \cite{Sim68} states that every closed, connected, embedded minimal submanifold $\Sigma^k\subset S^n$, $1\leq k<n$, that is not totally geodesic, satisfies the lower bound
\[
\Ind(\Sigma^k)\geq n+1\,.
\]
 The excluded totally geodesic $S^k$ has index $n-k$. The unstable directions arise from the normal projections of constant vectors in $\R^{n+1}$, or equivalently from ambient conformal variations. The natural next question is whether equality determines the submanifold.

For orientable minimal hypersurfaces in $S^n$, namely in the codimension one case, the stronger bound $n+2$ holds \cite[Theorem 2.2]{ElS93} and so equality in the previous estimate is impossible. For surfaces in $S^3$, Urbano \cite{Urb90}  characterized the equality case: a closed, connected, orientable minimal surface that is not totally geodesic has index at least five, with equality precisely for the Clifford torus $S^1(1/\sqrt2)\times S^1(1/\sqrt2)$, up to an ambient isometry. This classification became a decisive ingredient in Marques--Neves' proof of the Willmore conjecture \cite{MN14}. Their min--max argument bounds the index of a least-area minimal surface of positive genus by five; Urbano's theorem identifies it with the Clifford torus, whose area $2\pi^2$ supplies the sharp Willmore bound. Note that in codimension one, a global unit normal turns the Jacobi operator into a scalar operator. In $S^4$, such an operator acts on a rank-two normal bundle, and the scenario is a lot more complicated.

The present paper gives the corresponding equality classification in $S^4$: the Veronese projective plane is the unique closed embedded minimal surface attaining Simons' lower bound of five. In particular, the index alone determines its topology, without a prior genus or orientability assumption. As explained below, this classification also yields a sharp Willmore inequality and identifies the only possible links of strongly isolated singularities of three-dimensional minimal submanifolds in a generic five-dimensional ambient manifold.

Throughout, surfaces are smooth, closed, and connected, and the sphere has radius one. For a minimal immersion $F:\Sigma\to S^4\subset\R^5$, let $N\to\Sigma$ be its real rank-two normal bundle in $S^4$. Thus
\[
N_p=\{\nu\in T_{F(p)}S^4:\langle\nu,dF_p(X)\rangle=0\text{ for all }X\in T_p\Sigma\}.
\]
The second fundamental form is the symmetric bilinear bundle map $A:T\Sigma\times_\Sigma T\Sigma\to N$, i.\,e. the smooth family of symmetric bilinear maps $A:T_p\Sigma\times T_p\Sigma\to N_p$ for any $p\in\Sigma$. For any vector bundle $E$, we write $\Gamma(E)$ for its smooth sections. Then, the Jacobi operator $L:\Gamma(N)\to\Gamma(N)$ is defined as
\begin{equation}\label{eq:jacobi}
L:=\Delta^\perp+B+2I,
\end{equation}
where $B:N\to N$ is the fiberwise self-adjoint endomorphism $ B(V):=\sum_{i,j=1}^2\langle A(e_i,e_j),V\rangle A(e_i,e_j)$, for $(e_1,e_2)$ any local orthonormal tangent frame, and $\Delta^\perp:\Gamma(N)\to\Gamma(N)$ is the nonpositive normal connection Laplacian. The index $\Ind(\Sigma)$ counts the negative eigenvalues of $-L$, with real multiplicities on normal sections over $\Sigma$. The nullity is $\Null(\Sigma):=\dim_\R\ker L$. These quantities refer to the given immersion $F$, even when $F$ is not an embedding. 

Let $\Sym_0(\R^3)$ be the space of tracefree symmetric $3\times3$ real matrices, with inner product $\langle U,V\rangle=\tr(UV)$. Its unit sphere is an isometric copy of $S^4$. The Veronese projective plane $\Pi$ is the image of the embedding $F_0:\RP^2\to S^4\subset\Sym_0(\R^3)$ given by
\begin{equation}\label{eq:model}
[u]\longmapsto\sqrt{\frac32}\left(uu^T-\frac13 I\right),\qquad u\in S^2(1).
\end{equation}
\begin{theorem}[Index-five rigidity]\label{thm:main}
Let $\Sigma\subset S^4$ be a closed, connected, smooth, embedded minimal surface that is not totally geodesic. Then $\Ind(\Sigma)\geq5$, and equality holds if and only if $\Sigma=R\Pi$ for some $R\in O(5)$.
\end{theorem}
The totally geodesic case is a great $S^2$, which has index two; see Section~\ref{subsec:simons}. In particular, equality forces $\Sigma$ to be nonorientable and diffeomorphic to $\RP^2$. To our knowledge, this is the first complete classification of a (nonempty) equality case Simons' bound for fixed submanifold and ambient dimensions, among closed connected embedded minimal submanifolds, without additional topological or geometric hypotheses. 

The closest antecedent to our result is Karpukhin's characterization \cite{Kar21}, stated explicitly in \cite[Theorem 4.8]{KW24}: among minimal immersions of $\RP^2$ into $S^n$, $n\geq4$, having Morse index equal to $n+1$ characterizes the Veronese plane in a great $S^4$. Earlier contributions include Ejiri's estimate
\cite{Eji83}: a linearly full minimal immersion of
$S^2$ into $S^{2m}$ has Morse index at least
$2m(m+2)-6$. In $S^4$, this gives the sharp lower
bound ten, attained by the double cover of the
Veronese projective plane. The projective plane
itself has index five; the larger index of its
double cover reflects the availability of normal
variations that do not descend to $\RP^2$.
Ejiri's subsequent work \cite{Eji86} provides the
degree-three rigidity used in the projective-plane
classification: if the orientation double cover of a
minimal immersion of $\RP^2$ into $S^4$ has a horizontal
holomorphic twistor lift of degree three, the immersion
is congruent to the Veronese plane; see
\cite[proof of Theorem 4.8]{KW24}.
For a closed orientable superminimal surface in
$S^4$ that is not totally geodesic,
Montiel--Urbano \cite{MU97} obtained the formula
$\Ind(\Sigma)=\Area(\Sigma)/\pi+2g-2$,
where $g$ is the genus of $\Sigma$.
Since such a surface has area at least $12\pi$,
this gives $\Ind(\Sigma)\geq 2g+10$,
extending the lower bound ten from minimal
two-spheres to orientable superminimal surfaces
of arbitrary genus. Kusner--Wang \cite{KW24} proved the lower bound six for minimal tori in $S^4$, with equality precisely for a Clifford torus in a great $S^3$. Bernstein--Ketover \cite[Theorem 1.1]{BK26a} subsequently (and very recently) excluded equality $n+1$ for orientable codimension-two submanifolds in even-dimensional spheres.
There are also higher-dimensional equality examples: Rauchenberger--Semmelmann \cite{RS26} computed indices $8$, $14$, and $26$ for the standard Veronese embeddings of $\CP^2$, $\mathbb H\mathrm P^2$, and $\mathbb O\mathrm P^2$ into $S^7$, $S^{13}$, and $S^{25}$, respectively. 

In what follows, let us conveniently define the shifted operator $\Lt:\Gamma(N)\to\Gamma(N)$ by
\begin{equation}\label{eq:shifted}
\Lt:=L-2I=\Delta^\perp+B,
\end{equation}
and let us agree that the condition $\Lt\leq0$ means that $\langle\Lt V,V\rangle_{L^2}\leq0$ for every smooth real normal section $V$ over $\Sigma$. At a point $p\in \Sigma$, the curvature ellipse is the set
\[
\mathcal E_p:=\{A_p(v,v):v\in T_p\Sigma,\ |v|=1\}\subset N_p.
\]
We call it circular if it is a circle centered at zero, allowing radius zero. A minimal immersion into $S^4$ is called
\emph{superminimal} if its curvature ellipse is
circular at every point, allowing the degenerate
circle $\{0\}$ at points where $A=0$. This terminology was introduced by
Bryant \cite{Bry82}, although the study of these
surfaces predates his work. One key step towards to proof of our main result, Theorem \ref{thm:main} is the following slightly more technical statement:

\begin{theorem}[Superminimality from nonpositivity of the shifted operator]\label{thm:circular}
Let $F:\Sigma\to S^4$ be a smooth minimal immersion of a closed, connected surface. If $\Lt\leq0$, then either $F$ is superminimal, or $F(\Sigma)$ is contained in a great three-sphere. In particular, for an embedding, $\Lt\leq0$ implies that $F$ is superminimal.
\end{theorem}
In the nonorientable embedded case, the great-three-sphere alternative is excluded by topology. In the orientable case, $\Lt\leq0$ forces a surface contained in a great $S^3$ to be totally geodesic, hence superminimal as well. The description of superminimal surfaces by holomorphic curves is classical; see e.\,g. Eells--Salamon \cite{ES85}. Montiel--Urbano \cite{MU97} studied the relation between their index, area, and genus. Here we use their holomorphic normal-kernel formula to recover the topology from the index hypothesis.
\subsection{Consequences for Willmore energy and isolated singularities}
For a smooth immersion $F:\Sigma\to S^4$, define its Willmore energy by
\[
\mathcal W(F):=\int_\Sigma\left(1+\frac14|\mathbf H|^2\right)\,d\mu_\Sigma,
\]
where $\mathbf H:=\tr A$ is the mean curvature vector in $S^4$ and $d\mu_\Sigma$ is the induced area measure. Under stereographic projection from a pole outside $F(\Sigma)$, this energy becomes one quarter of the integral of the squared norm of the Euclidean mean curvature vector of the projected immersion, taken with respect to its induced Euclidean area measure. Kusner's conjecture \cite{Kus96} predicts that a stereographic image of the Veronese projective plane minimizes Willmore energy among closed surfaces in $\R^4$ not diffeomorphic to $S^2$. Combining Theorem~\ref{thm:main} with Bernstein--Ketover's minimal-filling theorem gives the following instance of that conjecture.

For an immersion into an oriented four-manifold, the \emph{Euler number of the normal bundle}, denoted by $e(N)$, is the signed number of zeros of a normal section transverse to the zero section. At each zero, choose tangent and normal orientations whose direct sum is the ambient orientation, and take the sign of the determinant of the section's derivative $T_p\Sigma\to N_p$. Reversing both orientations leaves this sign unchanged, so the definition also applies to nonorientable surfaces; see \cite[Definition 3.1 and Remark 3.2]{GR18}. This integer is distinct from $\chi(\Sigma)$.

\begin{corollary}[Sharp Willmore bound]\label{cor:willmore}
Let $F:\Sigma\to S^4$ be a smooth immersion of a closed, connected surface whose normal bundle has nonzero Euler number. Then $\mathcal W(F)\geq6\pi$, with equality if and only if $F$ is an embedding whose image is a conformal image of $\Pi$. Likewise, every smooth immersion of a closed, connected surface into $\R^4$ with nonzero normal Euler number has Euclidean Willmore energy at least $6\pi$, with equality if and only if it is an embedding whose image is a stereographic projection of a conformal image of $\Pi$, with projection pole outside that image.
\end{corollary}
\begin{proof}
An orientable embedding in $S^4$ has normal Euler number zero, since this number is its self-intersection number and $H_2(S^4;\Z)=0$. Bernstein--Ketover \cite[Theorem 1.3]{BK26b} and Theorem~\ref{thm:main} therefore give the lower bound $6\pi$ for embeddings with $e(N)\neq0$. Nonembedded immersions have Willmore energy at least $8\pi$ by Li--Yau \cite{LY82}.

For completeness, we justify the equality characterization with its conformal qualification. Suppose $\mathcal W(F)=6\pi$, so that $F$ is an embedding. Write
\[
\lambda_c(F):=\sup_{\psi\in\operatorname{Conf}(S^4)}\Area(\psi\circ F).
\]
Conformal invariance gives $\lambda_c(F)\leq\mathcal W(F)$, while \cite[Theorem 1.8]{BK26b}, together with Theorem~\ref{thm:main}, gives $\lambda_c(F)\geq6\pi$. Hence $\lambda_c(F)=6\pi$. For a smooth embedded surface, a conformal-area supremum strictly above $4\pi$ is attained; see also the compactness observation in the proof of \cite[Proposition 3.4]{BK26b}. Indeed, a sequence of conformal transformations escaping every compact set modulo isometries has limiting area at most $4\pi$, by concentration at a single point of the embedded surface. Choose $\psi$ attaining the supremum. Since
\[
\Area(\psi\circ F)=6\pi=\mathcal W(\psi\circ F),
\]
$\psi\circ F$ is minimal. It has least area among minimal embeddings with nonzero normal Euler number, because the Willmore lower bound applies to all of them. Thus \cite[Corollary 1.4]{BK26b} gives index five, and Theorem~\ref{thm:main} identifies its image with $\Pi$ up to an ambient isometry.

Conversely, $\Pi$ has area and Willmore energy $6\pi$. Conformal diffeomorphisms preserve Willmore energy and the absolute normal Euler number, so all its conformal images attain equality. Stereographic projection with pole outside the surface, and its inverse, preserve these quantities and embeddedness. This proves both the Euclidean bound and its equality characterization.
\end{proof}
No orientability or nonsphericity assumption is needed in Corollary~\ref{cor:willmore}: every orientable immersion with $e(N)\neq0$, including any such immersed sphere, is not embedded and thus satisfies the stronger bound $\mathcal W\geq8\pi$.

\begin{remark}[Scope of the normal Euler number hypothesis]
\label{rem:eulerscope}
For a closed nonorientable surface, let $k$ be its cross-cap number, namely the number of projective-plane summands, so that $\chi(\Sigma)=2-k$. By the Whitney--Massey theorem \cite[Section 2]{Mas69}, the possible normal Euler numbers of smooth embeddings into $S^4$ are exactly
\[
e(N)\in\{-2k+4j:j=0,1,\ldots,k\}.
\]
Thus the hypothesis $e(N)\neq0$ in Corollary~\ref{cor:willmore} is automatic for odd $k$, whereas for even $k$ both zero and nonzero values occur. In particular, an embedded $\RP^2$ ($k=1$), including the Veronese plane, has normal Euler number $\pm2$; an embedded Klein bottle ($k=2$) can have normal Euler number $-4$, $0$, or $4$. These examples concern smooth embeddings, without a minimality assumption.
\end{remark}

Our second application concerns the non-persistence theorem proven by the authors in \cite{CLW24}. A singular point of a stationary integral varifold is \emph{strongly isolated} if it has a multiplicity-one tangent cone whose link is smooth and embedded; in this setting the tangent cone is unique by Simon's theorem \cite{Simon83}. Write
\[
C(\Pi):=\{rx:r\geq0,\ x\in\Pi\}\subset\R^5.
\]
\begin{corollary}[Generic singularity models in dimension five]\label{cor:generic}
Let $M$ be a closed smooth five-manifold. For a Baire-generic smooth Riemannian metric $g$ on $M$, every $g$-stationary integral three-varifold with only strongly isolated singularities has, at each singular point $p$, a multiplicity-one tangent cone congruent to $C(\Pi)$ under a linear isometry $(T_pM,g_p)\simeq\R^5$.
\end{corollary}
\begin{proof}
Let $p$ be a singular point and let $\Lambda\subset S(T_pM,g_p)\simeq S^4$ be the smooth embedded minimal link of its multiplicity-one tangent cone. The link is connected: otherwise two of its components would be disjoint closed minimal surfaces in $S^4$, contrary to Frankel's intersection theorem \cite{Fra66}. Moreover, $\Lambda$ cannot be totally geodesic. Indeed, its cone would then be a multiplicity-one three-plane, and Allard's regularity theorem \cite{All72} would imply that $p$ is a regular point, a contradiction. Theorem 1.1 of \cite{CLW24}, with $n=3$ and $N=5$, forces $\Ind(\Lambda)=5$. Theorem~\ref{thm:main} therefore identifies $\Lambda$ with the Veronese projective plane up to an ambient isometry.
\end{proof}
As a result, the third alternative in the main theorem of \cite{CLW24} (namely, Theorem 1.1 therein) consists entirely of Veronese-cone singularities. This alternative cannot be removed: \cite[Corollary B]{CLW26} constructs minimal submanifolds modeled on the spherical suspension of $\Pi$ in $S^5$ for every sufficiently small smooth perturbation of the round metric. The present result shows that, in these dimensions, the persistent example realizes the only possible strongly isolated tangent-cone model for generic metrics. 

\subsection{Sketch of the proof of Theorem \ref{thm:main}}
There are four main steps. First, by Simons' work there holds the preliminary reduction (Proposition~\ref{prop:five}): for a surface that is not totally geodesic, index five implies
\[
\Lt\leq0,\qquad \ker\Lt=\{a^\perp:a\in\R^5\},\qquad \dim_\R\ker\Lt=5,
\]
where $a^\perp\in\Gamma(N)$ is the normal projection of the constant vector $a$.

Second, we construct a Dirac bundle over $\Sigma$ by pulling back
the spinor bundle of $S^4$ along $F$ and modifying its connection
according to the spinorial Gauss formula, as in
\eqref{eq:connection}; as a byproduct, this naturally determines an associated self-adjoint Dirac operator
$\mathscr D:\Gamma(\Sigma\times\C^4)\to
\Gamma(\Sigma\times\C^4)$.  We then derive a quadratic-form identity~\eqref{eq:trace} which relates every section $\Psi$ of this bundle to eight real normal fields $W_a[\Psi]$ of $\Sigma\subset S^4$ and, when $\Lt\leq0$, forces
\[
\spec\mathscr D\subset(-\infty,-1]\cup\{0\}\cup[1,\infty).
\]
Proposition~\ref{prop:paired} then shows that if $\ker\mathscr D=0$, then the surface lies in a great $S^3$. Its proof uses Riemann--Roch and Hodge theory to find a nonzero pair consisting of a harmonic tangent field $X$ and a normal field $w$ of $\Sigma\subset S^4$, satisfying $\mathcal A_\pm(X,w)=0$ for one of the operators defined in~\eqref{eq:augmented}, and thus produces suitable test sections satisfying the opposite inequality $\|\mathscr D\Psi\|^2_{L^2}\leq \|\Psi\|^2_{L^2}$. The resulting rigidity argument forces $w$ to be parallel, whence the conclusion. 
For a nonorientable embedded surface, this is topologically impossible; in the orientable case, $\Lt\leq0$ would force total geodesicity.

Third, applying the quadratic-form identity~\eqref{eq:trace} to a nonzero solution of $\mathscr D\Psi=0$ gives $\Qt(W_a[\Psi])=0$ for each of the eight associated normal fields. Corollary~\ref{cor:zeroshiftedform}, applied under $\Lt\leq0$, then gives $\Lt W_a[\Psi]=0$. These equations force superminimality (Proposition~\ref{prop:equality}), proving Theorem~\ref{thm:circular}. In the orientable case the normal complex structure then commutes with $\Lt$, making its real kernel even-dimensional. A five-dimensional kernel therefore forces nonorientability (Corollary~\ref{cor:orientableindex}).

Finally, the holomorphic description of the normal kernel, real descent to $\Sigma$, and Riemann--Roch express the real dimension of the kernel of the shifted normal Jacobi operator $\Lt$ as $2d+g-1$. Since $\dim_\R\ker\Lt=5$, this gives
\[
5=2d+g-1,\qquad d\geq3,
\]
where $g$ is the genus of the orientation double cover $\St$ and $d=\deg\Fh^*\mathcal O(1)$ is the degree of its horizontal holomorphic twistor lift $\Fh:\St\to\CP^3$. Thus $g=0$ and $d=3$. Proposition~\ref{prop:degree} therefore gives $\Sigma\cong\RP^2$, and Karpukhin's characterization \cite{Kar21}, stated in \cite[Theorem 4.8]{KW24}, identifies the surface with the Veronese plane. The same result supplies the converse $\Ind(\Pi)=5$.

The main contribution is the passage from nonpositivity of the shifted normal Jacobi operator to superminimality. In the proof, the quadratic-form identity~\eqref{eq:trace}, the paired augmented operators, and the equality argument in Sections~\ref{sec:clifford}--\ref{sec:circular} supply this passage. Their analytic ingredients are standard Dirac-bundle theory \cite{LM89}, the spinorial Gauss formula \cite{Bar98,BLR13}, Hodge theory, and Riemann--Roch. Once superminimality is established, we use the classical holomorphic normal-bundle calculation \cite{MU97} and twistor description \cite{Bry82,ES85}, with the deck symmetry kept explicit to descend to the nonorientable surface. Section~\ref{sec:degree} states these classical inputs and applies them without repeating their proofs; the remaining argument is the descent to the nonorientable surface and the determination of its topology.

\subsection{Organization}
Section~\ref{sec:prelim} fixes the geometric and analytic conventions and recalls Simons' normal fields and the index reduction in Section~\ref{subsec:simons}. Section~\ref{sec:clifford} constructs $\mathscr D$ and proves the quadratic-form identity~\eqref{eq:trace}, and Section~\ref{sec:circular} proves superminimality. Section~\ref{sec:degree} applies the classical superminimal formulas to determine the topology and completes the proof of Theorem~\ref{thm:main} using the known index characterization of minimal projective planes.

\section{Preliminaries}\label{sec:prelim}

\subsection{Curvature ellipses}
Equip $\Sigma$ with the induced metric $g_\Sigma:=F^*g_{S^4}$. The tangent and normal connections are here denoted by
\[
\nabla^\Sigma:\Gamma(T\Sigma)\to\Gamma(T^*\Sigma\otimes T\Sigma),\qquad
\nabla^\perp:\Gamma(N)\to\Gamma(T^*\Sigma\otimes N).
\]
At each $p\in\Sigma$ and for each normal vector $\nu\in N_p$, the \emph{shape operator in the direction $\nu$} is the self-adjoint linear map
\begin{equation}\label{eq:Weingarten}
A_\nu:T_p\Sigma\longrightarrow T_p\Sigma,\qquad
\langle A_\nu X,Y\rangle=\langle A_p(X,Y),\nu\rangle
\quad(X,Y\in T_p\Sigma).
\end{equation}
For a smooth normal field $\nu\in\Gamma(N)$, the same notation denotes the resulting bundle endomorphism $A_\nu:T\Sigma\to T\Sigma$.

The Euclidean connection $D^{\R^5}$ acts on sections of the trivial bundle $\Sigma\times\R^5$. For $X,Y\in\Gamma(T\Sigma)$ and $\nu\in\Gamma(N)$, the Gauss--Weingarten equations are
\begin{equation}\label{eq:gaussweingarten}
\begin{aligned}
D_X^{\R^5}dF(Y)&=dF(\nabla_X^\Sigma Y)+A(X,Y)-\langle X,Y\rangle F,\\
D_X^{\R^5}\nu&=-dF(A_\nu X)+\nabla_X^\perp\nu.
\end{aligned}
\end{equation}
We identify $T_p\Sigma$ with $dF_p(T_p\Sigma)$ when viewing tangent vectors as ambient vectors. Accordingly, $F^*TS^4=T\Sigma\oplus N$ under this identification.

For a local orthonormal tangent frame $(e_1,e_2)$ on $U\subset\Sigma$, write $A_{ij}:=A(e_i,e_j)\in\Gamma(N|_U)$. Minimality gives $A_{22}=-A_{11}$. At each point of $U$, the curvature ellipse is
\[
\{A_p(X,X):X\in T_p\Sigma,\ |X|=1\}=\{\cos(2t)A_{11}+\sin(2t)A_{12}:t\in\R\}.
\]
Fix an orientation of $S^4$. On an oriented patch of $\Sigma$, orient $N$ so that the tangent and normal orientations together give this ambient orientation. On such a patch $U\subset\Sigma$, let the orthogonal complex structures determined by the metric and these orientations be
\[
J_T:T\Sigma|_U\to T\Sigma|_U,\qquad J_N:N|_U\to N|_U.
\]
Thus $J_T^2=-\operatorname{id}_{T\Sigma}$ and $J_N^2=-\operatorname{id}_N$; in positively oriented orthonormal frames they satisfy $J_Te_1=e_2$ and $J_N\nu_1=\nu_2$. Reversing the tangent orientation also reverses the normal orientation, so both rotations change sign. 

For a connection on a bundle $E$, we use the curvature convention
\[
R(X,Y):=\nabla_X\nabla_Y-\nabla_Y\nabla_X-\nabla_{[X,Y]}\in\Gamma(\operatorname{End}(E)),
\qquad X,Y\in\Gamma(T\Sigma).
\]
For the normal bundle, we take $E=N$ and $\nabla=\nabla^\perp$, and denote its curvature by $R^\perp$. Thus $R^\perp(X,Y)\in\Gamma(\operatorname{End}(N))$ for tangent fields $X,Y$. With
\begin{equation}\label{eq:uv}
u:=A_{11},\qquad v:=A_{12},\qquad s:=|u|^2+|v|^2,\qquad\kappa:=2\det(u,v),
\end{equation}
where $\det(u,v)$ is computed in a positively oriented orthonormal normal frame, the Gauss and Ricci equations read
\begin{equation}\label{eq:gaussricci}
K=1-s,\qquad R^\perp(e_1,e_2)=-\kappa J_N.
\end{equation}
The first identity is the Gauss equation for the Gaussian curvature $K:\Sigma\to\R$; the second is the Ricci equation for the normal connection, with the curvature convention fixed above. Although $u$ and $v$ depend on the tangent frame, $s,\kappa:\Sigma\to\R$ are globally defined smooth functions. Indeed, simultaneous tangent and normal reflection leaves $2\det(A_{11},A_{12})$ unchanged.

\begin{lemma}[Criteria for a circular curvature ellipse]\label{lem:circularitycriteria}
At every point $p\in\Sigma$, the quantities in~\eqref{eq:uv} and the endomorphism $B_p:N_p\to N_p$ (defined right after \eqref{eq:jacobi}) satisfy
\begin{equation}\label{eq:curvaturebounds}
|A|^2=\tr B=2s,\qquad |\kappa|\leq s=\frac{|A|^2}{2}.
\end{equation}
Moreover, the following conditions are equivalent (including when $A_p=0$):
\begin{equation}\label{eq:circularitycriteria}
\begin{aligned}
\mathcal E_p\text{ is circular}
&\quad\Longleftrightarrow\quad |u|=|v|\ \text{ and }\ \langle u,v\rangle=0\\
&\quad\Longleftrightarrow\quad |\kappa|=s
\quad\Longleftrightarrow\quad B_p=s\operatorname{id}_{N_p}.
\end{aligned}
\end{equation}
\end{lemma}
\begin{proof}
Minimality and the definition of $B$ give
\[
|A|^2=2(|u|^2+|v|^2),\qquad
B_p\nu=2\langle u,\nu\rangle u+2\langle v,\nu\rangle v
\quad(\nu\in N_p),
\]
so $\tr B=2s$. Indeed, writing $w(t)=u\cos(2t)+v\sin(2t)$, we have
\[
|w(t)|^2
=\frac{|u|^2+|v|^2}{2}
+\frac{|u|^2-|v|^2}{2}\cos(4t)
+\langle u,v\rangle\sin(4t).
\]
Thus $|w(t)|$ is constant exactly when $|u|=|v|$ and
$\langle u,v\rangle=0$. Under these conditions, if the common
length is positive, $w$ parametrizes a circle centered at zero;
otherwise $u=v=0$, giving the circle of radius zero.

Also,
\[
s^2-\kappa^2=(|u|^2-|v|^2)^2+4\langle u,v\rangle^2,
\]
which proves the inequality and its equality criterion. Finally, in an orthonormal normal frame let $U$ be the $2\times2$ matrix with columns $u,v$. Then $B_p=2UU^T$. If $s>0$, the identity $B_p=sI$ is equivalent to $U^TU=(s/2)I$, again expressing orthogonality and equal lengths. If $s=0$, then $u=v=0$ and all the stated conditions hold.
\end{proof}

\subsection{Quadratic forms and elliptic operators.}\label{subsec:forms}
Let $d\mu_\Sigma$ be the area measure of $g_\Sigma$. This measure is defined without choosing an orientation. If $F$ is an embedding, it agrees, under identification by $F$, with the two-dimensional Hausdorff measure $d\mathscr H^2$ on $F(\Sigma)$. For a general immersion we integrate over the domain $\Sigma$; its pushforward therefore counts overlapping sheets with their multiplicities. The same convention applies to the orientation cover, equipped with the pullback metric. 

For a metric vector bundle $E\to\Sigma$, let $L^2(E)$ and $H^k(E)$ denote the spaces of square-integrable sections and Sobolev sections of order $k$, respectively. Their norms use $d\mu_\Sigma$ and the bundle metric; we abbreviate the $L^2$ norm by $\|\cdot\|_2$. The nonpositive normal connection Laplacian $\Delta^\perp:\Gamma(N)\to\Gamma(N)$ is
\[
\Delta^\perp V:=\sum_i\left(\nabla_{e_i}^\perp\nabla_{e_i}^\perp V-\nabla_{\nabla_{e_i}^{\Sigma}e_i}^\perp V\right).
\]
Keeping in mind the second variation formula for the area functional \cite{Sim68}, we introduce the following quadratic functionals $Q,\Qt:H^1(N)\to\R$:

\begin{equation}\label{eq:forms}
\begin{aligned}
Q(V)&:=\int_\Sigma\bigl(|\nabla^\perp V|^2-\langle BV,V\rangle-2|V|^2\bigr)\,d\mu_\Sigma,\\
\Qt(V)&:=\int_\Sigma\bigl(|\nabla^\perp V|^2-\langle BV,V\rangle\bigr)\,d\mu_\Sigma.
\end{aligned}
\end{equation}
By the definitions of $B$ and the shape operator,
\[
\langle BV,V\rangle
=\sum_{i,j=1}^2\langle A(e_i,e_j),V\rangle^2
=|A_V|^2.
\]
Thus, equivalently,
\[
Q(V)=\int_\Sigma
\bigl(|\nabla^\perp V|^2-|A_V|^2-2|V|^2\bigr)\,d\mu_\Sigma,
\qquad
\Qt(V)=\int_\Sigma
\bigl(|\nabla^\perp V|^2-|A_V|^2\bigr)\,d\mu_\Sigma.
\]
For $V\in\Gamma(N)$, integration by parts also gives
\[
Q(V)=\langle-LV,V\rangle_{L^2},\qquad
\Qt(V)=-\langle\Lt V,V\rangle_{L^2},
\]
where $L$ and $\Lt$ are defined in
\eqref{eq:jacobi} and \eqref{eq:shifted}, respectively.
Thus $\Lt\leq0$ is equivalent to $\Qt(V)\geq0$ for every $V\in H^1(N)$.

\medskip\noindent\textbf{Eigenvalues and spectrum.} The following facts are the standard elliptic spectral results on a closed manifold; see \cite[Chapter III, Section 5, especially Theorem 5.8]{LM89} and \cite[Sections 1.3--1.6]{Gil84}. The last two assertions use $\Lt=L-2I$ and $Q(V)=-\langle LV,V\rangle_{L^2}$. None of these statements requires orientability.

\begin{proposition}[Spectral properties of the normal operators]
\label{prop:normalspectrum}
The operators $L$ and $\Lt$ are self-adjoint on $L^2(N)$ with domain $H^2(N)$. Each has real eigenvalues of finite multiplicity, with no finite accumulation point, and admits a complete $L^2$-orthonormal basis of smooth real eigensections. There are no other spectral points. Moreover,
\[
\spec\Lt=\{\lambda-2:\lambda\in\spec L\},\qquad
\ker(\Lt-(\lambda-2)I)=\ker(L-\lambda I),
\]
so the shift preserves eigensections and multiplicities. The Morse index is
\[
\Ind(\Sigma)=\sum_{\lambda\in\spec L,\ \lambda>0}
\dim_\R\ker(L-\lambda I).
\]
\end{proposition}

\begin{corollary}[Vanishing of the shifted quadratic form]
\label{cor:zeroshiftedform}
Suppose $\Lt\leq0$. For every $V\in\Gamma(N)$,
\[
\Qt(V)=0\quad\Longleftrightarrow\quad\Lt V=0.
\]
\end{corollary}
\begin{proof}
Expanding $V$ in the eigenbasis of $\Lt$ gives $\Qt(V)=\sum_j(-\lambda_j)|a_j|^2$, with $\lambda_j\leq0$, so the form vanishes exactly when all components with $\lambda_j<0$ vanish.
\end{proof}
This corollary will turn equality in the quadratic-form identity~\eqref{eq:trace} into individual equations for normal fields.

\medskip\noindent\textbf{Principal symbols and ellipticity.}
Let $E,G\to\Sigma$ be smooth real or complex vector
bundles, and let $P:\Gamma(E)\to\Gamma(G)$ be a linear
differential operator of order one.
For $p\in\Sigma$, $\xi\in T_p^*\Sigma$, and $v\in E_p$,
choose a real smooth function $f$ with $df_p=\xi$
and a smooth section $u$ of $E$ with $u(p)=v$.
The \emph{principal symbol} of $P$ is the fiber map
$\sigma_1(P)(p,\xi):E_p\to G_p$ defined by
\begin{equation}\label{eq:symboldefinition}
\sigma_1(P)(p,\xi)v
:=\bigl(P(fu)-fP(u)\bigr)(p).
\end{equation}
This expression is independent of the choices of
$f$ and $u$. In local coordinates and bundle frames, if $P$ has order $m$ and local expression
\[
P=\sum_{|\alpha|\leq m}a_\alpha(x)\partial^\alpha,
\]
its principal symbol is
\[
\sigma_m(P)(p,\xi)
:=\sum_{|\alpha|=m}a_\alpha(p)\xi^\alpha:E_p\to G_p.
\]
Here $\alpha=(\alpha_1,\alpha_2)$ is a multi-index,
$|\alpha|=\alpha_1+\alpha_2$,
$\partial^\alpha=\partial_1^{\alpha_1}
\partial_2^{\alpha_2}$, and
$\xi^\alpha=\xi_1^{\alpha_1}\xi_2^{\alpha_2}$.
These local expressions define an intrinsic fiber map.
Our convention omits the factors of $i$ used in
the Fourier convention.
The operator $P$ is called \emph{elliptic} if
$\sigma_m(P)(p,\xi)$ is invertible for every
$p\in\Sigma$ and every nonzero $\xi\in T_p^*\Sigma$.

\medskip\noindent\textbf{Fredholm operators and their indices.}
A bounded linear map $P:U\to V$ between real or complex Hilbert spaces is \emph{Fredholm} if its range is closed and both $\ker P$ and $\operatorname{coker}P:=V/\operatorname{ran}P$ are finite-dimensional. Over the scalar field $\mathbb F=\R$ or $\C$, its index is
\[
\ind_{\mathbb F}P:=\dim_{\mathbb F}\ker P-\dim_{\mathbb F}\operatorname{coker}P.
\]
Thus a \emph{complex index} counts complex dimensions. For a $\C$-linear Fredholm map viewed as $\R$-linear, $\ind_\R P=2\ind_\C P$. An elliptic differential operator of order $m$ between bundles over a closed surface is Fredholm as a map $H^m\to L^2$. Indeed, elliptic estimates and compact Sobolev inclusion give a finite-dimensional kernel and closed range; applying the same argument to the formal adjoint gives a finite-dimensional cokernel, identified with its smooth adjoint kernel. These facts do not require orientability; see \cite[Sections 1.3--1.6]{Gil84}.

\medskip
\noindent\textbf{Harmonic one-forms, Hodge theory and nonorientable surfaces.} A smooth real one-form $\alpha$ is \emph{harmonic} if $(dd^*+d^*d)\alpha=0$, equivalently $d\alpha=d^*\alpha=0$ on the closed surface. Here $d^*$ denotes the formal adjoint with respect to the induced metric and area measure. Real harmonic one-forms represent $H^1(\Sigma;\R)$, including when $\Sigma$ is nonorientable \cite[Section 1.5, Lemma 1.5.3]{Gil84}. Writing $b_j(\Sigma):=\dim_\R H^j(\Sigma;\R)$ for its Betti numbers, a closed connected nonorientable surface has $b_0=1$ and $b_2=0$, hence
\begin{equation}\label{eq:betti}
b_1(\Sigma)=1-\chi(\Sigma).
\end{equation}
When $\Sigma$ is nonorientable, let $p:\St\to\Sigma$ be the orientation double cover, and let $\tau:\St\to\St$ be its nontrivial deck transformation. Since $\Sigma$ is connected and nonorientable, $\St$ is connected. When using this cover, $g$ denotes the genus of the oriented surface $\St$, not of $\Sigma$. The identity
$2\chi(\Sigma)=\chi(\St)=2-2g$
then gives $\chi(\Sigma)=1-g$. If $k$ denotes instead the nonorientable genus, namely the number of cross-caps, then $\chi(\Sigma)=2-k$ and $k=g+1$. In particular, $\RP^2$ has cross-cap number one and orientation-cover genus zero, whereas the Klein bottle has cross-cap number two and orientation-cover genus one. For an orientable surface of genus $h$, we use $\chi(\Sigma)=2-2h$ and $b_1(\Sigma)=2h=2-\chi(\Sigma)$.

\subsection{Holomorphic line bundles and real descent}
An oriented Riemannian surface has isothermal coordinates and hence a Riemann surface structure. We apply this to $\St$ with its pullback metric. For a Hermitian complex line bundle $E\to\St$, a unitary connection
$\nabla^E:\Gamma(E)\longrightarrow\Omega^1(\St;E)$
has a $(0,1)$ part
$\bar\partial_E:\Gamma(E)\longrightarrow\Omega^{0,1}(\St;E)$.
Here $\Omega^{a,b}(\St;E)$ denotes smooth $E$-valued forms of type $(a,b)$. Since complex dimension is one, the integrability condition is automatic, and $\bar\partial_E$ defines a holomorphic structure on $E$. Its kernel is the complex vector space $H^0(E)$ of holomorphic sections. The curvature $F_E\in\Omega^2(\St;i\R)$ is an imaginary-valued two-form, and our degree convention is
\begin{equation}\label{eq:degreeconvention}
\deg E=\frac{i}{2\pi}\int_{\St}F_E.
\end{equation}
Write $h^0(E):=\dim_\C H^0(E)$. The Riemann--Roch theorem, in its operator-index form, gives
\begin{equation}\label{eq:rr}
\ind_\C\bar\partial_E=\deg E+1-g.
\end{equation}
Here the index is taken for the $\C$-linear Fredholm map
$\bar\partial_E:H^1(E)\longrightarrow L^2(\Lambda^{0,1}T^*\St\otimes E)$,
where $\Lambda^{0,1}T^*\St$ is the bundle of complex covectors of type $(0,1)$; see \cite[Theorem 28]{Bob11}.

 Recall that a map $T:V\to W$ between complex vector spaces is \emph{anti-linear} if
\[
T(az+bw)=\bar a\,Tz+\bar b\,Tw
\qquad(a,b\in\C,\ z,w\in V).
\]
For an anti-linear involution $T:V\to V$ (i.\,e. an anti-linear map satisfying $T^2=I$) we define 
$V^T:=\{v\in V:T(v)=v\}\,$
and note this is a real vector subspace of $V$, generally not a complex subspace.

\medskip
\noindent\textbf{Real descent from the orientation cover.}
Let $\pi:\widetilde\Sigma\to\Sigma$ be the orientation
double cover of a closed nonorientable surface, with
deck transformation $\tau:\widetilde\Sigma\to\widetilde\Sigma$.
Suppose that a complex vector bundle $E\to\widetilde\Sigma$
is equipped with smooth anti-linear isomorphisms
\[
\rho_p:E_p\longrightarrow E_{\tau(p)},
\qquad
\rho_{\tau(p)}\circ\rho_p=I.
\]
Anti-linearity here means that
$\rho_p(zv)=\overline z\,\rho_p(v)$ for $z\in\C$ and $v\in E_p$.

These data define a real vector bundle $E_\R\to\Sigma$.
Indeed, for $\pi^{-1}(x)=\{p,\tau(p)\}$, set
\[
(E_\R)_x
:=
\left\{
(v_p,v_{\tau(p)})\in E_p\oplus E_{\tau(p)}
:\ v_{\tau(p)}=\rho_p(v_p)
\right\}.
\]
The condition is unchanged upon exchanging the two lifts.
Projection onto either component is a $\R$-linear
isomorphism with the corresponding fiber of $E$;
thus $\operatorname{rank}_\R E_\R
=2\operatorname{rank}_\C E$.
Locally, choosing one sheet of the cover identifies
$E_\R$ with $E$ on that sheet, regarded as a real bundle.
The smooth maps $\rho_p$ identify the two choices and
give $E_\R$ its smooth bundle structure.

Consequently, sections of $E_\R$ correspond to sections
$s\in\Gamma(E)$ satisfying
\begin{equation}\label{eq:realdescentsections}
s(\tau(p))=\rho_p(s(p)).
\end{equation}
Equivalently, they are the fixed sections of the
anti-linear involution
\[
T:\Gamma(E)\longrightarrow\Gamma(E),
\qquad
(Ts)(p):=\rho_{\tau(p)}(s(\tau(p))).
\]
We call such sections \emph{invariant}.
The same correspondence holds for Sobolev sections.

\begin{lemma}[Real descent]\label{lem:real}
For $j=0,1$, let $E_j\to\widetilde\Sigma$ carry maps
$\rho_j$ as above, and denote by $E_{j,\R}\to\Sigma$
and $T_j$ the associated real bundles and involutions
on sections. Let
$P:\Gamma(E_0)\longrightarrow\Gamma(E_1)$
be a $\C$-linear elliptic differential operator
of order $m$ satisfying $PT_0=T_1P$.
Then $P$ preserves invariant sections and induces
a $\R$-linear elliptic operator
$P_\R:\Gamma(E_{0,\R})\longrightarrow\Gamma(E_{1,\R})$.
For the Fredholm maps
\[
P:H^m(E_0)\longrightarrow L^2(E_1),
\qquad
P_\R:H^m(E_{0,\R})\longrightarrow L^2(E_{1,\R}),
\]
one has
\[
\dim_\R\ker P_\R=\dim_\C\ker P,
\qquad
\dim_\R\operatorname{coker}P_\R
=\dim_\C\operatorname{coker}P.
\]
In particular,
\[
\operatorname{ind}_\R P_\R
=\operatorname{ind}_\C P.
\]
\end{lemma}

\begin{proof}
The identity $PT_0=T_1P$ ensures that $P$ maps invariant
sections to invariant sections. 
For $j=0,1$, define the $\R$-linear isomorphism
\[
I_j:\Gamma(E_{j,\R})\longrightarrow
\Gamma(E_j)^{T_j},
\qquad
(I_j u)(p):=\operatorname{pr}_p\bigl(u(\pi(p))\bigr),
\]
where $\operatorname{pr}_p:(E_{j,\R})_{\pi(p)}
\to(E_j)_p$ selects the component indexed by $p$,
and
$\Gamma(E_j)^{T_j} :=\{s\in\Gamma(E_j):T_js=s\}$.
Define
\[
P_\R:=I_1^{-1}\circ P\circ I_0.
\]
Under the local bundle identifications obtained by
choosing a sheet of the cover, $P_\R$ has the same
differential expression as $P$, acting on the underlying
real bundles. Its principal symbol is therefore
invertible at every nonzero real cotangent vector,
because that of $P$ is. Thus $P_\R$ is elliptic, and both displayed Sobolev
maps are Fredholm.

If a complex vector space $V$ carries an anti-linear
involution $T$, denote its fixed (real) subspace by
$V^T:=\{v\in V:Tv=v\}$.
Every $v\in V$ has the unique decomposition
\[
v=\frac{v+Tv}{2}
+i\,\frac{v-Tv}{2i},
\]
with both displayed components in $V^T$.
Thus, when $V$ is finite-dimensional,
\[
\dim_\R V^T=\dim_\C V.
\]
Applying this observation to $\ker P$, whose fixed
subspace corresponds to $\ker P_\R$, proves the
kernel formula.

The involution $T_1$ also induces an anti-linear
involution $[y]\mapsto[T_1y]$ on $\operatorname{coker}P$.
Every fixed class has an invariant representative,
namely $(y+T_1y)/2$. Moreover, if an invariant $y$
belongs to the range of $P$, say $y=Px$, then
\[
y=P\left(\frac{x+T_0x}{2}\right),
\]
so it has an invariant preimage.
Hence $\operatorname{coker}P_\R$ identifies with
the fixed subspace of $\operatorname{coker}P$.
The same dimension identity proves the cokernel
formula and therefore the index formula.
\end{proof}

\subsection{The holomorphic quartic}
\begin{lemma}\label{lem:quartic}
Let $F:\Sigma\to S^4$ be a minimal immersion of a connected oriented surface, and let $z$ be a local conformal coordinate. Set $A_{zz}:=A(\partial_z,\partial_z)$, that is a locally defined section of $N\otimes_\R\C$. The quartic differential
\[\langle A_{zz},A_{zz}\rangle\,dz^4\]
is a globally defined holomorphic section of $\Kc_\Sigma^{\otimes4}$, where $\Kc_\Sigma=(T^{1,0}\Sigma)^*$ and the normal inner product is extended complex bilinearly. Its zeros are exactly the circular points. It either vanishes identically or has isolated zeros.
\end{lemma}
\begin{proof}
In a conformal coordinate $z$, minimality gives $A_{z\bar z}=0$ and the tangent connection satisfies $\nabla_{\bar z}^\Sigma\partial_z=\nabla_z^\Sigma\partial_{\bar z}=0$. The Codazzi equation therefore implies $\nabla_{\bar z}^\perp A_{zz}=0$. Metric compatibility gives holomorphicity. In an oriented orthonormal frame, the coefficient is a nonzero real factor times
$|u|^2-|v|^2-2i\langle u,v\rangle.$
The zero criterion follows from the definition of the ellipse. The final assertion follows from Lemma~\ref{lem:circularitycriteria} and well-known facts on holomorphic functions.
\end{proof}

\subsection{Simons' normal fields and the index reduction}
\label{subsec:simons}

We recall the normal-field construction underlying Simons'
index estimate \cite[Section 5.1]{Sim68}.
For $a\in\R^5$, define $V_a\in\Gamma(N)$ by
$V_a(p):=\operatorname{proj}_{N_p}a$. In our conventions,
Simons' normal-field identity \cite[Lemma 5.1.3]{Sim68} reads $-LV_a=-2V_a$ namely
\begin{equation}\label{eq:canonicaljacobi}
\Lt V_a=0.
\end{equation}
For a minimal surface that is not totally geodesic, these
fields depend injectively on $a$. We record the resulting
index reduction in the form needed below.

\begin{proposition}[Simons' index estimate and its equality consequence]
\label{prop:five}
Let $F:\Sigma\to S^4$ be a closed, connected minimal immersion
that is not totally geodesic. The linear map
\[
\R^5\longrightarrow\ker\Lt,\qquad a\longmapsto V_a,
\]
is injective. Consequently, $\Ind(\Sigma)\geq5$.
If $\Ind(\Sigma)=5$, then $\ker\Lt=\{V_a:a\in\R^5\}$ and
$\Lt\leq0$.
\end{proposition}

\begin{remark}\label{rem:totallygeodesic}
A closed, connected totally geodesic immersion into $S^4$
has image a great $S^2$ and defines a covering of that sphere.
Since $S^2$ is simply connected, this covering has one sheet
\cite[Theorem 1.38]{Hat02}. In particular, $\Sigma$ is
diffeomorphic to $S^2$, so no closed nonorientable surface
admits such an immersion. The great $S^2\subset S^4$ has
Morse index two, with its two negative directions given by
the constant normal fields.
\end{remark}

\section{An ambient Clifford connection and a quadratic-form identity}\label{sec:clifford}
Throughout this section, $F:\Sigma\to S^4$ is a smooth minimal immersion of a closed, connected surface, equipped with its induced metric. No orientability assumption is made. Fix the ambient orientation as in Section~\ref{sec:prelim}. No assumption on $\Lt$ is made until the spectral conclusion below.

The underlying construction we are about to present belongs to the standard spin geometry of submanifolds. B\"ar \cite[Sections 1--2]{Bar98} describes the ambient spinor bundle and the spinorial Gauss formula; Bayard--Lawn--Roth \cite[Sections 2.2--2.3]{BLR13} give the corresponding formulas for surfaces in four-dimensional space forms. The Dirac operator and its Weitzenb\"ock formula are standard \cite[Chapter II, Sections 5 and 8]{LM89}; the calculation specific to our argument is their relation to the shifted normal quadratic form $\Qt$ in Proposition~\ref{prop:trace}.

\subsection{Clifford matrices, the connection, and the Dirac operator}
We identify $T_p\Sigma$ with $dF_p(T_p\Sigma)\subset\R^5$. Then
\[
T_{F(p)}S^4=dF_p(T_p\Sigma)\oplus N_p=F(p)^\perp
\qquad(p\in\Sigma).
\]
We shall associate to each vector
$v\in T_{F(p)}S^4$ a matrix $c_p(v)$, depending
$\R$-linearly on $v$, such that
\[
c_p(v)^2=-|v|^2I
\qquad(v\in T_{F(p)}S^4).
\]
Applying this identity to $v+w$ shows that matrices
associated to orthogonal vectors must anticommute.
Thus an orthonormal basis of $T_{F(p)}S^4$ requires
four matrices whose squares are $-I$ and which
anticommute in pairs. Such matrices exist on $\C^4$,
and no smaller complex vector space admits them
\cite[Chapter I, Section 5]{LM89}.
We therefore take
\[
\SSb:=\Sigma\times\C^4\longrightarrow\Sigma,
\]
equipped with its standard Hermitian metric. A section of $\SSb$ is identified with a smooth $\C^4$-valued function on $\Sigma$. The formulas below use the immersion and fixed matrices on $\C^4$, and require no choice of a spin or Pin structure on $\Sigma$.

Choose five Hermitian matrices $\gamma_1,\ldots,\gamma_5\in\operatorname{End}_\C(\C^4)$ satisfying
\[
\gamma_\alpha\gamma_\beta+\gamma_\beta\gamma_\alpha
=2\delta_{\alpha\beta}I\quad(1\leq\alpha,\beta\leq5),
\qquad \gamma_1\gamma_2\gamma_3\gamma_4\gamma_5=I,
\]
where $I$ is the identity on $\C^4$. Such matrices exist by the standard Clifford matrix construction \cite[Chapter I, Section 5]{LM89}; changing the sign of one matrix, if necessary, gives the displayed product condition. For $a=(a_1,\ldots,a_5)\in\R^5$, define the $\R$-linear map
\[
\gamma:\R^5\longrightarrow\operatorname{End}_\C(\C^4),
\qquad \gamma(a):=\sum_{\alpha=1}^5a_\alpha\gamma_\alpha.
\]
The matrix relations imply
\[
\gamma(u)\gamma(v)+\gamma(v)\gamma(u)
=2\langle u,v\rangle I\qquad(u,v\in\R^5).
\]
For each $p\in\Sigma$, define
\[
c_p:T_{F(p)}S^4\longrightarrow\operatorname{End}_\C(\SSb_p),
\qquad c_p(v):=\gamma(F(p))\gamma(v).
\]
Thus $c_p$ acts on both tangent vectors $dF_p(X)$, with $X\in T_p\Sigma$, and normal vectors $\nu\in N_p$. It is $\R$-linear and satisfies
\[
c_p(v)^*=-c_p(v),\qquad
c_p(v)c_p(w)+c_p(w)c_p(v)=-2\langle v,w\rangle I
\quad(v,w\in T_{F(p)}S^4).
\]
Here $*$ denotes the Hermitian adjoint. These identities follow from $v,w\perp F(p)$ and $\gamma(F(p))^2=I$. To improve readability, we suppress $p$ when applying $c$ to vector fields along $F$; in particular, $c(Z)$ means $c(dF(Z))$ for a tangent field $Z$.

We next define the connection. For $Z\in\Gamma(T\Sigma)$ and $\Psi\in\Gamma(\SSb)$, let $d_Z\Psi$ denote the ordinary directional derivative of the $\C^4$-valued function $\Psi$ along $Z$. Choose a local orthonormal tangent frame $(e_1,e_2)$, and recall that $A(Z,e_j)\in\Gamma(N)$ for $j=1,2$. Define $S_Z\in\Gamma(\operatorname{End}_\C(\SSb))$ and a connection $\nabla:\Gamma(\SSb)\to\Gamma(T^*\Sigma\otimes\SSb)$ by
\begin{equation}\label{eq:connection}
\begin{aligned}
S_Z&:=\frac12c(Z)-\frac12\sum_{j=1}^2c(e_j)c(A(Z,e_j)),\\
\nabla_Z\Psi&:=d_Z\Psi+S_Z\Psi.
\end{aligned}
\end{equation}
The sum defining $S_Z$ is independent of the orthonormal frame, as verified in Lemma~\ref{lem:cliffordconnection} below.

To express $S_Z$ in terms of $\gamma$ rather than $c$,
observe that
\[
c_p(v)c_p(w)=-\gamma(v)\gamma(w)
\qquad(p\in\Sigma,\ v,w\in T_{F(p)}S^4).
\]
Substituting this identity and $c(Z)=\gamma(F)\gamma(dF(Z))$ into \eqref{eq:connection} gives
\[
S_Z=\frac12\gamma(F)\gamma(dF(Z))
+\frac12\sum_{j=1}^2\gamma(dF(e_j))\gamma(A(Z,e_j))
\qquad(Z\in\Gamma(T\Sigma)).
\]
Equivalently, the connection is obtained from the spherical spin connection $\widetilde\nabla_Z\Psi=d_Z\Psi+\tfrac12c(Z)\Psi$ by subtracting $\tfrac12\sum_{j=1}^2c(e_j)c(A(Z,e_j))\Psi$. This is the spinorial Gauss formula \cite[Section 2, equation (9)]{Bar98}; see also \cite[Section 2.3]{BLR13}. Strictly speaking, both papers derive such a formula by assuming that the ambient manifold and the submanifold admits a spin structure. However, we do not need that assumption, and we just define the connection this way.

Write $\nabla^{T\Sigma\oplus N}:=\nabla^\Sigma\oplus\nabla^\perp$. Compatibility of $\nabla$ with this connection means that differentiating $c(v)\Psi$ obeys the Leibniz rule in \eqref{eq:compatibility}, with $v$ differentiated using $\nabla^{T\Sigma\oplus N}$ and $\Psi$ using $\nabla$.
In expressions involving $\gamma$, tangent vectors continue to be identified with their images under $dF$.

\begin{lemma}[Clifford compatibility]\label{lem:cliffordconnection}
Formula~\eqref{eq:connection} defines a global Hermitian connection on $\SSb$. For $Z\in\Gamma(T\Sigma)$, $v\in\Gamma(T\Sigma\oplus N)$, and $\Psi\in\Gamma(\SSb)$, it satisfies
\begin{equation}\label{eq:compatibility}
\nabla\gamma(F)=0,\qquad
\nabla_Z(c(v)\Psi)=c(\nabla_Z^{T\Sigma\oplus N}v)\Psi+c(v)\nabla_Z\Psi.
\end{equation}
In particular, restriction of $c$ to $T\Sigma$ makes $(\SSb,c,\nabla)$ a Dirac bundle over $\Sigma$.
\end{lemma}
\begin{proof}
The sum defining $S_Z$ is an orthogonal contraction of the tensor $A$, so it is unchanged under any orthonormal tangent-frame change, including reflections. Each summand in $S_Z$ is skew-Hermitian: $c(Z)$ is skew-Hermitian, and the Clifford products of perpendicular tangent and normal vectors are skew-Hermitian. Thus the connection is global and preserves the Hermitian metric.
The spinorial Gauss formula cited above gives compatibility with the direct-sum connection. In our matrix conventions, for
$Z\in\Gamma(T\Sigma)$ and
$v\in\Gamma(T\Sigma\oplus N)$, the identities to check are
\[
d_Z\gamma(F)+[S_Z,\gamma(F)]=0,\qquad
d_Z\gamma(v)+[S_Z,\gamma(v)]
=\gamma(\nabla_Z^{T\Sigma\oplus N}v),
\]
where brackets denote matrix commutators. They follow by substituting the displayed formula for $S_Z$ and the Gauss--Weingarten equations~\eqref{eq:gaussweingarten}; the radial and second-fundamental-form terms cancel. Applying the Leibniz rule to $c(v)=\gamma(F)\gamma(v)$ gives~\eqref{eq:compatibility}.
\end{proof}

For a local orthonormal tangent frame $(e_1,e_2)$, write $\nabla_i:=\nabla_{e_i}$ for $i=1,2$. The associated \emph{Dirac operator} is the first-order operator $\mathscr D:\Gamma(\SSb)\to\Gamma(\SSb)$ defined by
\begin{equation}\label{eq:dirac}
\mathscr D:=\sum_i c(e_i)\nabla_i.
\end{equation}
This is precisely the standard Clifford contraction of the covariant derivative; see \cite[Chapter II, Section 5, equation (5.0)]{LM89}. It depends on the specified Clifford action and connection, and is not being defined as the intrinsic spin Dirac operator of a chosen spin structure on $\Sigma$.

Define the \emph{ambient chirality operator} by
\[
\Omega(p):=\gamma(F(p))
\in\operatorname{End}_\C(\SSb_p)
\qquad(p\in\Sigma).
\]

\begin{lemma}[Algebraic properties of chirality]
\label{lem:chiralityalgebra}
The endomorphism $\Omega$ is a Hermitian involution:
\[
\Omega^*=\Omega,\qquad \Omega^2=I.
\]
For every $p\in\Sigma$ and $v\in T_{F(p)}S^4$,
\[
\Omega(p)c_p(v)=-c_p(v)\Omega(p).
\]
If $(e_1,e_2)$ and $(\nu_1,\nu_2)$ are local
orthonormal tangent and normal frames such that
$(F,e_1,e_2,\nu_1,\nu_2)$ is positively oriented in
$\R^5$, then
\[
\Omega=c(e_1)c(e_2)c(\nu_1)c(\nu_2).
\]

The orthogonal projections
\[
P_\pm:=\tfrac12(I\pm\Omega)
\]
have constant complex rank two. Their images
\[
\SSb^\pm:=\operatorname{im}P_\pm=\ker(\Omega\mp I)
\]
are smooth complex subbundles, called the
\emph{chirality subbundles}, and
$\SSb=\SSb^+\oplus\SSb^-$ orthogonally.
In a suitable smooth local unitary frame,
\[
\Omega=\begin{pmatrix}I_2&0\\0&-I_2\end{pmatrix}.
\]
\end{lemma}

\begin{proof}
Since the matrices $\gamma_\alpha$ are Hermitian and
$|F|=1$, we have
\[
\Omega^*=\gamma(F)^*=\gamma(F)=\Omega,
\qquad
\Omega^2=\gamma(F)^2=I.
\]
Fix $p\in\Sigma$ and $v\in T_{F(p)}S^4$.
The orthogonality $v\perp F(p)$ gives
\[
\gamma(F(p))\gamma(v)=-\gamma(v)\gamma(F(p)).
\]
Using $c_p(v)=\gamma(F(p))\gamma(v)$, we obtain
\[
\Omega(p)c_p(v)=\gamma(v),
\qquad
c_p(v)\Omega(p)=-\gamma(v),
\]
which proves the anticommutation identity.

To establish the frame formula, define the alternating
multilinear map
\[
\mathcal A(b_1,\ldots,b_5)
:=
\frac1{5!}\sum_{\sigma\in S_5}
\operatorname{sgn}(\sigma)\,
\gamma(b_{\sigma(1)})\cdots\gamma(b_{\sigma(5)})
\qquad(b_1,\ldots,b_5\in\R^5),
\]
where $S_5$ denotes the permutation group.
For an orthonormal basis, the matrices $\gamma(b_j)$
anticommute in pairs, so
\[
\mathcal A(b_1,\ldots,b_5)
=\gamma(b_1)\cdots\gamma(b_5).
\]
Since $\mathcal A$ is alternating and multilinear,
its value on a positively oriented orthonormal basis
equals its value on the standard basis, namely
$\gamma_1\cdots\gamma_5=I$.
Applying this to $(F,e_1,e_2,\nu_1,\nu_2)$ yields
\[
\gamma(F)\gamma(e_1)\gamma(e_2)
\gamma(\nu_1)\gamma(\nu_2)=I.
\]
Moreover, $c(v)c(w)=-\gamma(v)\gamma(w)$ whenever
$v,w\in T_{F(p)}S^4$. Therefore
\[
\begin{aligned}
c(e_1)c(e_2)c(\nu_1)c(\nu_2)
&=\gamma(e_1)\gamma(e_2)\gamma(\nu_1)\gamma(\nu_2)\\
&=\gamma(F)=\Omega.
\end{aligned}
\]
This formula is independent of the chosen frames
with the specified combined orientation. In particular,
simultaneously reversing the tangent and normal
orientations leaves the product unchanged.

Because $\Omega$ is a Hermitian involution, the
endomorphisms $P_\pm=\tfrac12(I\pm\Omega)$ are
orthogonal projections onto its $\pm1$ eigenspaces,
and $P_++P_-=I$.
At each point, choose a unit vector $v\in T_{F(p)}S^4$.
The anticommutation identity shows that $c_p(v)$
exchanges these eigenspaces. Since $c_p(v)^2=-I$,
this map is invertible. The two eigenspaces therefore
have equal complex dimension, necessarily two.

Finally, $P_\pm$ are smooth projections of constant
rank, so their images are smooth rank-two complex
subbundles. Choosing smooth local frames of these
images and applying Gram--Schmidt gives orthonormal
frames of $\SSb^+$ and $\SSb^-$.
Their union is a smooth local unitary frame of $\SSb$
in which $\Omega=\operatorname{diag}(I_2,-I_2)$.
\end{proof}

The following proposition collects the standard analytic facts we will use. They apply because Lemma~\ref{lem:cliffordconnection} makes $\SSb$ a Dirac bundle: ellipticity and formal self-adjointness follow from \cite[Chapter II, Lemma 5.1 and Proposition 5.3]{LM89}, and the domain and spectral assertions follow from elliptic regularity and the spectral theorem \cite[Chapter III, Section 5, especially Theorem 5.8]{LM89}. No orientability assumption is required.

\begin{proposition}[Standard analytic properties of the Dirac operator]
\label{prop:diracanalytic}
The operator $\mathscr D$ is elliptic and formally self-adjoint. Its closure on the complex Hilbert space $L^2(\SSb)$ is self-adjoint, with domain
$\mathcal D(\mathscr D)=H^1(\SSb)$.
Its spectrum consists of real eigenvalues of finite multiplicity, with no finite accumulation point, and its smooth eigensections form an orthonormal basis of $L^2(\SSb)$.
\end{proposition}

Every spectral statement below refers to this operator. In particular, $\lambda$ is an eigenvalue if $\mathscr D\Psi=\lambda\Psi$ for some nonzero $\Psi\in H^1(\SSb)$; as in Section~\ref{subsec:forms}, its spectrum consists precisely of its eigenvalues. Passing to the underlying real Hilbert space doubles each eigenspace dimension and leaves the spectrum unchanged.

\begin{lemma}[Chirality and the curvature term]\label{lem:diracidentities}
The ambient chirality operator $\Omega$ is parallel
with respect to $\nabla$. Moreover, $\Omega\mathscr D=-\mathscr D\Omega$ whence $\mathscr D$ maps
$\Gamma(\SSb^\pm)$ into $\Gamma(\SSb^\mp)$.
The Weitzenb\"ock identity is
\begin{equation}\label{eq:weitzenbock}
\mathscr D^2=\nabla^*\nabla+\frac K2I-\frac\kappa2\Omega.
\end{equation}
Here $\nabla^*$ is the formal $L^2$ adjoint of the
connection $\nabla$ on $\SSb$ defined in
\eqref{eq:connection}. 
\end{lemma}
\begin{proof}
Parallelness follows from $\Omega=\gamma(F)$ and
\eqref{eq:compatibility}. For every
$\Psi\in\Gamma(\SSb)$, Lemma~\ref{lem:chiralityalgebra}
therefore gives
\[
\mathscr D(\Omega\Psi)
=\sum_{i=1}^2 c(e_i)\Omega\nabla_{e_i}\Psi
=-\Omega\mathscr D\Psi.
\]
Consequently, $\mathscr D$ maps sections of either
chirality subbundle to sections of the other.
For $X,Y\in\Gamma(T\Sigma)$, let
$R^{\SSb}(X,Y)\in\Gamma(\operatorname{End}_\C(\SSb))$
denote the curvature of $\nabla$, with the convention
of Section~\ref{sec:prelim}.
Choose local orthonormal tangent and normal frames
$(e_1,e_2)$ and $(\nu_1,\nu_2)$ such that
$(F,e_1,e_2,\nu_1,\nu_2)$ is positively oriented in $\R^5$.
We claim that
\[
R^{\SSb}(e_1,e_2)
=-\frac K2c(e_1)c(e_2)
-\frac\kappa2c(\nu_1)c(\nu_2).
\]

Set $E:=T\Sigma\oplus N$, equipped with the connection
$\nabla^E:=\nabla^\Sigma\oplus\nabla^\perp$, and denote
its curvature by $R^E$. Differentiating the compatibility
identity \eqref{eq:compatibility} twice gives
\[
[R^{\SSb}(X,Y),c(v)]=c(R^E(X,Y)v)
\qquad
(X,Y\in\Gamma(T\Sigma),\ v\in\Gamma(E)),
\]
where $[U,V]:=UV-VU$ for endomorphisms $U,V$.
Thus the curvature of $\nabla$ acts on Clifford
multiplication by the curvature of $\nabla^E$.

For the chosen local orthonormal frames, the Gauss
and Ricci equations \eqref{eq:gaussricci} give
\[
\begin{aligned}
R^E(e_1,e_2)e_1&=-Ke_2,
& R^E(e_1,e_2)e_2&=Ke_1,\\
R^E(e_1,e_2)\nu_1&=-\kappa\nu_2,
& R^E(e_1,e_2)\nu_2&=\kappa\nu_1.
\end{aligned}
\]
Define
\[
C:=-\frac K2c(e_1)c(e_2)
   -\frac\kappa2c(\nu_1)c(\nu_2).
\]
The Clifford relations imply the pointwise identity
\[
[c(a)c(b),c(v)]
=2\langle a,v\rangle c(b)-2\langle b,v\rangle c(a)
\qquad(a,b,v\in E_p,\ p\in\Sigma).
\]
Applying it to the four frame vectors yields
\[
\begin{aligned}
{}[C,c(e_1)]&=-Kc(e_2),
& [C,c(e_2)]&=Kc(e_1),\\
[C,c(\nu_1)]&=-\kappa c(\nu_2),
& [C,c(\nu_2)]&=\kappa c(\nu_1).
\end{aligned}
\]
Consequently,
$D:=R^{\SSb}(e_1,e_2)-C$
commutes with each of these four Clifford matrices.

The standard matrix description of the complex
Clifford algebra in dimension four
\cite[Chapter I, Section 5]{LM89} says that these
four matrices and their products span
$\operatorname{End}_\C(\C^4)$.
Hence $D$ commutes with every complex $4\times4$
matrix and must be a scalar multiple of $I$.

Finally, both $R^{\SSb}(e_1,e_2)$ and $C$ have zero
complex trace. Indeed, if matrices $U,V$ anticommute,
cyclicity gives
$\tr(UV)=\tr(VU)=-\tr(UV)$,
and hence $\tr(UV)=0$.
The expression for $S_Z$ in terms of $\gamma$ is a sum
of such products, since $F\perp dF(Z)$ and
$dF(e_j)\perp A(Z,e_j)$. Thus $\tr S_Z=0$.
In the product trivialization,
\[
R^{\SSb}(X,Y)
=d_XS_Y-d_YS_X-S_{[X,Y]}+[S_X,S_Y]
\qquad(X,Y\in\Gamma(T\Sigma)),
\]
where $d_XS_Y$ denotes differentiation of the
matrix-valued function $S_Y$ along $X$.
Taking traces shows that $R^{\SSb}(X,Y)$ is tracefree.
The same cyclicity argument gives $\tr C=0$.
Therefore $D$ is scalar and tracefree, hence zero,
proving the asserted curvature formula.

The general Dirac-bundle Weitzenb\"ock formula
\cite[Chapter II, Theorem 8.2]{LM89} gives
$
\mathscr D^2
=\nabla^*\nabla+c(e_1)c(e_2)R^{\SSb}(e_1,e_2),
$ but,
since $(c(e_1)c(e_2))^2=-I$ and
$c(e_1)c(e_2)c(\nu_1)c(\nu_2)=\Omega$, we obtain
\[
c(e_1)c(e_2)R^{\SSb}(e_1,e_2)
=\frac K2I-\frac\kappa2\Omega,
\]
which proves \eqref{eq:weitzenbock}.
Both $\kappa$ and $\Omega$ are unchanged under
simultaneous reversal of the tangent and normal
orientations, so the identity holds globally,
including when $\Sigma$ is nonorientable.
\end{proof}

A smooth section $\Psi\in\Gamma(\SSb)$ will be called \emph{Dirac-harmonic} if $\mathscr D\Psi=0$. This terminology refers specifically to this linear Dirac equation; it does not mean that $\Psi$ is parallel or lies in the kernel of the connection Laplacian. Harmonic one-forms, used separately below, retain the Hodge-theoretic meaning given in Section~\ref{subsec:forms}.

\subsection{The quadratic-form identity}
Let $\varepsilon_1,\ldots,\varepsilon_4$ be the standard
coordinate basis of $\C^4$. Identifying sections of
$\SSb=\Sigma\times\C^4$ with $\C^4$-valued functions,
define
\[
\phi_a(p):=\varepsilon_a,\qquad
\phi_{a+4}(p):=i\varepsilon_a
\qquad(p\in\Sigma,\ a=1,\ldots,4).
\]
On $\C^4$, we use the real inner product
\[
\langle z,w\rangle
:=\operatorname{Re}\sum_{k=1}^4z_k\overline{w_k}
\qquad(z,w\in\C^4).
\]
Then $(\phi_1,\ldots,\phi_8)$ is an orthonormal frame
of the underlying real vector bundle.
All eight sections satisfy $d_Z\phi_a=0$ for every
$Z\in\Gamma(T\Sigma)$, since their coordinates are
constant. Their covariant derivatives with respect
to $\nabla$ are computed in the following lemma.

\begin{lemma}[The fixed product frame]\label{lem:constantsections}
For $a\in\{1,\ldots,8\}$ and $i\in\{1,2\}$, set $\nabla_i:=\nabla_{e_i}$, $S_i:=S_{e_i}$. Then
\begin{equation}\label{eq:constantsections}
\nabla_i\phi_a=S_i\phi_a,\qquad \mathscr D\phi_a=-\phi_a.
\end{equation}
\end{lemma}
\begin{proof}
The first identity follows from $d\phi_a=0$
and~\eqref{eq:connection}. Contracting the definition
of $S_i$ with $c(e_i)$ gives the second, since
symmetry and minimality imply
$\sum_{i,j}c(e_i)c(e_j)c(A_{ij})=0$
and $\sum_i c(e_i)c(e_i)=-2I$.
\end{proof}

 For each $a\in\{1,\ldots,8\}$, define the $\R$-linear map $W_a:\Gamma(\SSb)\to\Gamma(N)$ by the pointwise pairing
\begin{equation}\label{eq:fields}
\langle W_a[\Psi](p),\nu\rangle=\langle c_p(\nu)\Psi(p),\phi_a(p)\rangle\qquad(p\in\Sigma,\ \nu\in N_p).
\end{equation}
The Riesz representation theorem in $N_p$ gives a unique vector $W_a[\Psi](p)$, and these vectors vary smoothly with $p$. Thus the construction produces normal sections on $\Sigma$ itself, independently of any choice of normal frame. \begin{proposition}[Quadratic-form identity and spectral gap]\label{prop:trace} Let $\tilde{Q}$ be the quadratic form defined by \eqref{eq:forms}.
For every $\Psi\in\Gamma(\SSb)$, the normal fields defined in~\eqref{eq:fields} satisfy
\begin{equation}\label{eq:trace}
\sum_{a=1}^8\Qt(W_a[\Psi])=2\bigl(\|\mathscr D\Psi\|_2^2-\langle\mathscr D\Psi,\Psi\rangle_{L^2}\bigr).
\end{equation}
If, in addition, $\Lt\leq0$, then
\begin{equation}\label{eq:gap}
\spec\mathscr D\subset(-\infty,-1]\cup\{0\}\cup[1,\infty).
\end{equation}
\end{proposition}
\begin{proof}
\emph{The quadratic-form identity.} 
At a fixed point, choose an orthonormal normal frame
$(\nu_1,\nu_2)$. Since $c(\nu_\alpha)^*c(\nu_\alpha)=I$ for
$\alpha=1,2$, multiplication by either normal frame
vector preserves the norm. By \eqref{eq:fields} and orthonormality
of $(\phi_1,\ldots,\phi_8)$,
\begin{equation}\label{eq:parseval}
\begin{aligned}
\sum_{a=1}^8|W_a[\Psi]|^2
&=\sum_{\alpha=1}^2\sum_{a=1}^8
  \langle c(\nu_\alpha)\Psi,\phi_a\rangle^2\\
&=\sum_{\alpha=1}^2|c(\nu_\alpha)\Psi|^2
=2|\Psi|^2,\\[4pt]
\sum_{a=1}^8\langle BW_a[\Psi],W_a[\Psi]\rangle
&=\sum_{i,j=1}^2\sum_{a=1}^8
  \langle A_{ij},W_a[\Psi]\rangle^2\\
&=\sum_{i,j=1}^2\sum_{a=1}^8
  \langle c(A_{ij})\Psi,\phi_a\rangle^2\\
&=\sum_{i,j=1}^2|c(A_{ij})\Psi|^2
=|A|^2|\Psi|^2
=2s|\Psi|^2.
\end{aligned}
\end{equation}
Here $A_{ij}:=A(e_i,e_j)$ for an orthonormal tangent
frame $(e_1,e_2)$ and we used
$|A|^2=2s$ from \eqref{eq:curvaturebounds}.

Fix $p\in\Sigma$. Choose local orthonormal tangent
and normal frames satisfying
\[
\nabla^\Sigma_{e_i}e_j(p)=0,\qquad
\nabla^\perp_{e_i}\nu_\alpha(p)=0
\qquad(i,j,\alpha\in\{1,2\}),
\]
with $(F,e_1,e_2,\nu_1,\nu_2)$ positively oriented.
All pointwise calculations below are at $p$.
Write $\nabla_i:=\nabla_{e_i}$ and
$S_i:=S_{e_i}$.

To begin with, we differentiate the defining identity
\[
\langle W_a[\Psi],\nu_\alpha\rangle
=\langle c(\nu_\alpha)\Psi,\phi_a\rangle.
\]
By \eqref{eq:compatibility} and our choice of normal
frame,
\[
\nabla_i(c(\nu_\alpha)\Psi)(p)
=c(\nu_\alpha)\nabla_i\Psi(p).
\]
Moreover, Lemma~\ref{lem:cliffordconnection} states
that $\nabla$ preserves the Hermitian metric.
Taking real parts gives
\[
e_i\langle\Psi,\Phi\rangle
=\langle\nabla_i\Psi,\Phi\rangle
+\langle\Psi,\nabla_i\Phi\rangle
\qquad(\Psi,\Phi\in\Gamma(\SSb)).
\]
Therefore, at $p$,
\[
\begin{aligned}
\langle\nabla_i^\perp W_a[\Psi],\nu_\alpha\rangle
&=e_i\langle c(\nu_\alpha)\Psi,\phi_a\rangle\\
&=\langle c(\nu_\alpha)\nabla_i\Psi,\phi_a\rangle
 +\langle c(\nu_\alpha)\Psi,\nabla_i\phi_a\rangle\\
&=\langle c(\nu_\alpha)\nabla_i\Psi,\phi_a\rangle
 +\langle c(\nu_\alpha)\Psi,S_i\phi_a\rangle\\
&=\langle c(\nu_\alpha)\nabla_i\Psi
       -S_i c(\nu_\alpha)\Psi,\phi_a\rangle.
\end{aligned}
\]
The last equality uses $S_i^*=-S_i$.
Summing squared components in the orthonormal frame
$(\phi_a)_{a=1}^8$ therefore yields
\[
\sum_{a=1}^8|\nabla^\perp W_a[\Psi]|^2
=\sum_{i,\alpha=1}^2
\bigl|c(\nu_\alpha)\nabla_i\Psi
      -S_i c(\nu_\alpha)\Psi\bigr|^2.
\]
We compute the coefficients in this square explicitly.
For brevity, within this calculation, let us set
\[
E_i:=c(e_i),\qquad N_\alpha:=c(\nu_\alpha),
\qquad U:=c(u),\qquad V:=c(v),
\]
where $u:=A_{11}$ and $v:=A_{12}$.
Minimality gives $A_{22}=-u$, so
\[
S_1=\tfrac12(E_1-E_1U-E_2V),\qquad
S_2=\tfrac12(E_2-E_1V+E_2U).
\]

First, the Clifford relations give
\[
\sum_{\alpha=1}^2N_\alpha^*E_iN_\alpha=-2E_i,
\qquad
\sum_{\alpha=1}^2N_\alpha^*E_jc(n)N_\alpha=0
\quad(n\in N_p).
\]
Indeed, each $N_\alpha$ anticommutes with $E_i$
and satisfies $N_\alpha^2=-I$, proving the first
identity. For the second, anticommutation with $E_j$
reduces the sum to
$E_j\sum_\alpha N_\alpha c(n)N_\alpha$.
Writing $n=n_1\nu_1+n_2\nu_2$ shows that this sum
vanishes. Consequently,
\[
\sum_{\alpha=1}^2N_\alpha^*S_iN_\alpha=-E_i
\qquad(i=1,2).
\]

Next, write $u=u_1\nu_1+u_2\nu_2$ and
$v=v_1\nu_1+v_2\nu_2$. By (the definition)
$\kappa=2(u_1v_2-u_2v_1)$,
$UV-VU=\kappa N_1N_2$.
Using $U^2=-|u|^2I$, $V^2=-|v|^2I$, and
$\Omega=E_1E_2N_1N_2$, we first obtain
\[
(E_1U+E_2V)^2=(E_1V-E_2U)^2
=-sI-E_1E_2(UV-VU)=-sI-\kappa\Omega.
\]
Expanding the displayed expressions for $S_1$ and $S_2$
therefore gives
\[
\begin{aligned}
4S_1^2&=-(1+s)I-\kappa\Omega-2E_1E_2V,\\
4S_2^2&=-(1+s)I-\kappa\Omega+2E_1E_2V.
\end{aligned}
\]
Since $S_i^*=-S_i$, adding these identities yields
\[
\sum_{i=1}^2S_i^*S_i
=\frac{1+s}{2}I+\frac\kappa2\Omega.
\]
Finally, $N_\alpha^*N_\alpha=I$ and
$\Omega N_\alpha=-N_\alpha\Omega$, so
\[
\begin{aligned}
\sum_{i,\alpha=1}^2N_\alpha^*S_i^*S_iN_\alpha
&=\sum_{\alpha=1}^2
 N_\alpha^*\left(\frac{1+s}{2}I
                +\frac\kappa2\Omega\right)N_\alpha\\
&=(1+s)I-\kappa\Omega.
\end{aligned}
\]

We now expand
\[
\sum_{i,\alpha=1}^2
\bigl|N_\alpha\nabla_i\Psi-S_iN_\alpha\Psi\bigr|^2,
\]
which equals $\sum_{a=1}^8|\nabla^\perp W_a[\Psi]|^2$
by the identity established above. The first terms
sum to $2|\nabla\Psi|^2$, since each $N_\alpha$
preserves the norm. The cross terms satisfy
\[
\begin{aligned}
\sum_{i,\alpha=1}^2
\langle N_\alpha\nabla_i\Psi,S_iN_\alpha\Psi\rangle
&=\sum_{i=1}^2
 \left\langle\nabla_i\Psi,
 \left(\sum_{\alpha=1}^2N_\alpha^*S_iN_\alpha\right)
 \Psi\right\rangle\\
&=-\sum_{i=1}^2\langle\nabla_i\Psi,E_i\Psi\rangle\\
&=\langle\mathscr D\Psi,\Psi\rangle.
\end{aligned}
\]
Together with the preceding computation of the
quadratic terms, this gives
\[
\sum_{a=1}^8|\nabla^\perp W_a[\Psi]|^2
=2|\nabla\Psi|^2+(1+s)|\Psi|^2
-\kappa\langle\Omega\Psi,\Psi\rangle
-2\langle\mathscr D\Psi,\Psi\rangle.
\]
Since $p$ was arbitrary, this identity holds throughout
$\Sigma$.

Subtracting the potential term as per \eqref{eq:parseval}
and using $K=1-s$, we obtain
\[
\begin{aligned}
\sum_{a=1}^8\Qt(W_a[\Psi])
=\int_\Sigma\Bigl(
2|\nabla\Psi|^2+K|\Psi|^2
-\kappa\langle\Omega\Psi,\Psi\rangle
-2\langle\mathscr D\Psi,\Psi\rangle
\Bigr)\,d\mu_\Sigma.
\end{aligned}
\]
On the other hand, integrating the Weitzenb\"ock
identity \eqref{eq:weitzenbock} gives
\[
2\|\mathscr D\Psi\|_{L^2}^2
=\int_\Sigma\Bigl(
2|\nabla\Psi|^2+K|\Psi|^2
-\kappa\langle\Omega\Psi,\Psi\rangle
\Bigr)\,d\mu_\Sigma.
\]
Here we used formal self-adjointness of $\mathscr D$
and the definition of $\nabla^*$.
Combining these two identities proves
\eqref{eq:trace}.

\medskip\noindent\emph{Consequence for the spectrum.} The left-hand side of \eqref{eq:trace} is a sum of nonnegative quantities when $\Lt\leq0$. Let $\Psi\in\Gamma(\SSb)$ be a nonzero eigensection with $\mathscr D\Psi=\lambda\Psi$. Substituting this section into the quadratic-form identity~\eqref{eq:trace} gives
\[
0\leq\sum_{a=1}^8\Qt(W_a[\Psi])
=2\lambda(\lambda-1)\|\Psi\|_{L^2}^2.
\]
Since $\Psi\neq0$, we obtain $\lambda(\lambda-1)\geq0$. Since $\Omega$ pairs $\lambda$ and $-\lambda$, applying the same inequality to $-\lambda$ excludes both intervals $(-1,0)$ and $(0,1)$. This proves~\eqref{eq:gap}.
\end{proof}

\section{From nonpositivity of the shifted operator to superminimality}\label{sec:circular}

In this section we prove Theorem~\ref{thm:circular}.
The proof uses tangent vector fields on $\Sigma$
whose metric-dual one-forms are harmonic, following
a strategy employed in second-variation estimates
by Ros \cite{Ros06} and Urbano \cite[Section 2]{Urb13}.

\subsection{Geometric rigidity when the Dirac operator is injective}
\begin{proposition}\label{prop:paired}
Let $F:\Sigma\to S^4$ be a closed, connected minimal immersion with $\Lt\leq0$. Then either $\ker \mathscr D\neq0$, or there is a nonzero parallel normal field $w$
such that $A_w:T\Sigma\to T\Sigma$ vanishes (see \eqref{eq:Weingarten}).
 In the latter case, $w$ is a fixed vector in $\R^5$ and $F(\Sigma)$ lies in the great sphere $w^\perp\cap S^4$.
\end{proposition}

\begin{proof}
Suppose that $\ker\mathscr D=0$. We will prove the second alternative.
By \eqref{eq:gap},
\begin{equation}\label{eq:Dsquare}
\|\mathscr D\Psi\|_2^2-\|\Psi\|_2^2\geq 0 \qquad \forall\, \Psi\in H^1(\SSb) \,.
\end{equation}
We seek geometric test sections for which the sum of these differences
is nonpositive.

\medskip\noindent\emph{Step 1: A general test section.}
Let $X\in\Gamma(T\Sigma)$ and $w\in\Gamma(N)$ be arbitrary smooth fields,
and write $X^\flat:=g_\Sigma(X,\cdot)$ for the metric-dual one-form.
Recall the projections $P_\pm$ from Lemma~\ref{lem:chiralityalgebra}
and the constant real orthonormal frame $(\phi_a)_{a=1}^8$ of $\SSb$.
For either choice of sign, consider
\[
\Psi_a:=c(X+w)P_\pm\phi_a\in\Gamma(\SSb),
\qquad a=1,\ldots,8.
\]
We first claim that
\begin{equation}\label{eq:generaltest}
\begin{aligned}
\mathscr D\Psi_a & = \biggl[-(\diver X)I + (dX^\flat)(e_1,e_2)c(e_1)c(e_2)\\
 &\qquad\qquad + (c(X)-2S_X) + \sum_i c(e_i)c(\nabla_i^\perp w)\biggr]P_\pm\phi_a 
 + c(w)P_\mp\phi_a ,
\end{aligned}
\end{equation}
where $S_X$ is defined in \eqref{eq:connection}.

To prove this, recall that $P_\pm$ are parallel, $c(e_i)$ interchanges
the two chirality subbundles, and
$S_i-\tfrac12c(e_i)=-\tfrac12\sum_jc(e_j)c(A(e_i,e_j))$
preserves them. Thus \eqref{eq:constantsections} gives
\[
\nabla_i(P_\pm\phi_a)
=\tfrac12c(e_i)P_\mp\phi_a
 +\bigl(S_i-\tfrac12c(e_i)\bigr)P_\pm\phi_a.
\]
Clifford compatibility \eqref{eq:compatibility} now yields
\begin{equation}\label{eq:testexpansion}
\begin{aligned}
\mathscr D\Psi_a = \sum_i c(e_i)c(\nabla^\Sigma_{e_i} X+\nabla_i^\perp w)P_\pm\phi_a & +\frac12\sum_i c(e_i)c(X+w)c(e_i)P_\mp\phi_a\\
 & + \sum_i c(e_i)c(X+w)
       \bigl(S_i-\tfrac12c(e_i)\bigr)P_\pm\phi_a.
\end{aligned}
\end{equation}
In tangent dimension two, the Clifford relations imply
\[
\sum_i c(e_i)c(X)c(e_i)=0,
\qquad
\sum_i c(e_i)c(w)c(e_i)=2c(w).
\]
Symmetry of $A$ and minimality also give
\[
\sum_i c(e_i)\bigl(S_i-\tfrac12c(e_i)\bigr)
=-\frac12\sum_{i,j}c(e_i)c(e_j)c(A(e_i,e_j))
=\frac12c(\tr A)=0.
\]
Consequently,
\begin{align*}
\sum_i c(e_i)c(X)\bigl(S_i-\tfrac12c(e_i)\bigr)
&=-2\sum_i\langle X,e_i\rangle
       \bigl(S_i-\tfrac12c(e_i)\bigr)
=c(X)-2S_X,\\
\sum_i c(e_i)c(w)\bigl(S_i-\tfrac12c(e_i)\bigr)&=0.
\end{align*}
Finally, expanding the tangent derivative gives
\begin{align*}
\sum_i c(e_i)c(\nabla^\Sigma_{e_i} X)
& = -\sum_i\langle\nabla^\Sigma_{e_i} X,e_i\rangle I 
  + \bigl(\langle\nabla^\Sigma_{e_1} X,e_2\rangle
        -\langle\nabla^\Sigma_{e_2} X,e_1\rangle\bigr)c(e_1)c(e_2) \\
& = -(\diver X)I+(dX^\flat)(e_1,e_2)c(e_1)c(e_2).
\end{align*}
Substituting these identities into \eqref{eq:testexpansion} proves
\eqref{eq:generaltest}. Although written in local frames, the contraction
is independent of those frames and requires no orientation of $\Sigma$.

\medskip\noindent\emph{Step 2: Sum of $|\mathscr D \Psi_a|^2-|\Psi_a|^2$.}
We now take $X^\flat$ to be harmonic, while keeping $w$ arbitrary.
Then $dX^\flat=0$ and $d^*X^\flat=-\diver X=0$, so the first two terms
in \eqref{eq:generaltest} vanish and
\[
\mathscr D\Psi_a
=\biggl(c(X)-2S_X+\sum_i c(e_i)c(\nabla_i^\perp w)\biggr)P_\pm\phi_a
 +c(w)P_\mp\phi_a.
\]
For an endomorphism field $U$, use the real Hilbert--Schmidt norms
\[
|U|_{\mathrm{HS},\R}^2:=\sum_{a=1}^8|U\phi_a|^2,
\qquad
\|U\|_{\mathrm{HS},2}^2
:=\int_\Sigma |U|_{\mathrm{HS},\R}^2\,d\mu_\Sigma.
\]
We claim that the sum of the energy differences is
\begin{equation}\label{eq:normpair}
\begin{aligned}
\sum_{a=1}^8\bigl(\|\mathscr D\Psi_a\|_2^2-\|\Psi_a\|_2^2\bigr)
= \left\|\biggl(c(X)-2S_X+\sum_i c(e_i)c(\nabla_i^\perp w)\biggr)
P_\pm\right\|_{\mathrm{HS},2}^2-4\|X\|_2^2 \,.
\end{aligned}
\end{equation}
Indeed, the two endomorphisms in the formula for $\mathscr D\Psi_a$
vanish on $\SSb^\mp$ and $\SSb^\pm$, respectively. Their real
Hilbert--Schmidt inner product is therefore zero, as one sees in a real
orthonormal basis adapted to $\SSb=\SSb^+\oplus\SSb^-$.
This cancellation occurs after summing over $a$.
Since $P_\pm$ have real rank four and $c(v)^*c(v)=|v|^2I$, we obtain
\begin{align*}
\sum_{a=1}^8|\mathscr D\Psi_a|^2
&=\left|\biggl(c(X)-2S_X+\sum_i c(e_i)c(\nabla_i^\perp w)\biggr)
P_\pm\right|_{\mathrm{HS},\R}^2+4|w|^2,\\
\sum_{a=1}^8|\Psi_a|^2&=4(|X|^2+|w|^2).
\end{align*}
Subtracting and integrating proves \eqref{eq:normpair}.

\medskip\noindent\emph{Step 3.}
We seek a nonzero pair $(X,w)$, with $X^\flat$ harmonic, such that
\begin{equation}\label{eq:normalcorrection}
\biggl(c(X)-2S_X+\sum_i c(e_i)c(\nabla_i^\perp w)\biggr)P_\pm=0
\end{equation}
for one choice of sign. Together with \eqref{eq:Dsquare} and
\eqref{eq:normpair}, this equation will force $X=0$ in the next step.

We formulate \eqref{eq:normalcorrection} as a Fredholm problem: Let
\[
\HH:=\{X\in\Gamma(T\Sigma):dX^\flat=0,\ d^*X^\flat=0\}
\]
be the real space of harmonic tangent fields, and define the real
subbundles $\mathcal B_\pm\subset\operatorname{End}_\C(\SSb)$ by
\[
(\mathcal B_\pm)_p
:=\Span_\R\{c_p(X)c_p(\nu)P_\pm(p):X\in T_p\Sigma,\ \nu\in N_p\}.
\]
Define
\begin{equation}\label{eq:Cpm}
\begin{aligned}
C_\pm: \Gamma(N)\longrightarrow\Gamma(\mathcal B_\pm), \qquad 
C_\pm w := \sum_{i=1}^2c(e_i)c(\nabla_i^\perp w)P_\pm.
\end{aligned}
\end{equation}
This contraction is independent of the orthonormal tangent frame.
We will use its Sobolev extension to define
\begin{equation}\label{eq:augmented}
\begin{aligned}
\mathcal A_\pm&:\HH\oplus H^1(N)\longrightarrow L^2(\mathcal B_\pm),\\
\mathcal A_\pm(X,w)&:=\bigl(c(X)-2S_X\bigr)P_\pm+C_\pm w.
\end{aligned}
\end{equation}
Here $c(X)-2S_X=\sum_i c(e_i)c(A(X,e_i))$, so its product with
$P_\pm$ is a section of $\mathcal B_\pm$.
For smooth $X\in\HH$ and $w\in\Gamma(N)$, the formula in Step~2 becomes
\begin{equation}\label{eq:Dpair}
\mathscr D\bigl(c(X+w)P_\pm\phi_a\bigr)
=\mathcal A_\pm(X,w)\phi_a+c(w)P_\mp\phi_a.
\end{equation}
Thus the required pairs are precisely the nonzero elements of
$\ker\mathcal A_\pm$.

First note that $C_\pm$ are elliptic. In fact, choose orthonormal tangent
and normal frames with the prescribed combined orientation. The Clifford
relations give
\[
c(e_2)c(\nu_1)P_\pm=\mp c(e_1)c(\nu_2)P_\pm,
\qquad
c(e_2)c(\nu_2)P_\pm=\pm c(e_1)c(\nu_1)P_\pm.
\]
Hence $\mathcal B_\pm$ have real rank two, with local frame
$c(e_1)c(\nu_1)P_\pm,c(e_1)c(\nu_2)P_\pm$.
These endomorphisms are independent because $c(e_1)c(w)$ is invertible
for every nonzero $w\in N_p$ and cannot annihilate $P_\pm(p)\neq0$.
The same relations imply
\begin{equation}\label{eq:Clocal}
\begin{aligned}
C_+w&=c(e_1)c(\nabla_1^\perp w-J_N\nabla_2^\perp w)P_+,\\
C_-w&=c(e_1)c(\nabla_1^\perp w+J_N\nabla_2^\perp w)P_-.
\end{aligned}
\end{equation}
Here $J_N$ is the normal rotation from Section~\ref{sec:prelim}.
For $0\neq\xi\in T_p^*\Sigma$, the real principal symbol, with the
convention \eqref{eq:symboldefinition}, is
\[
N_p\longrightarrow(\mathcal B_\pm)_p,
\qquad w\longmapsto c_p(\xi^\#)c_p(w)P_\pm(p).
\]
It is injective by the same invertibility argument and is therefore an
isomorphism between real two-dimensional fibers. The Fredholm theorem
recalled in Section~\ref{subsec:forms} now applies to
$C_\pm:H^1(N)\to L^2(\mathcal B_\pm)$.

We compute both indices so that the normal-degree terms cancel.
First suppose $\Sigma$ is nonorientable. On its orientation cover $\St$,
of genus $g$, let $E:=(p^*N,J_N)$ be the normal complex line and $d_E:=\deg E$.

We describe explicitly how the lifted operators
$C_\pm$ correspond to the two complex components
of the normal connection.
Choose a positively oriented orthonormal tangent
frame $(e_1,e_2)$ on an open subset of $\St$, and
let $(\theta^1,\theta^2)$ be its dual coframe.
Multiplication by $i$ on $E$ means applying $J_N$. For this choice of tangent frame, the maps
\[
w\longmapsto c(e_1)c(w)P_\pm
\]
identify $E$ with $p^*\mathcal B_\pm$ as real
vector bundles over this open set.
These identifications depend on the tangent frame.
The following maps, however, are independent of
that choice and therefore define global bundle
isomorphisms:
\[
\begin{aligned}
\Phi_+:\ p^*\mathcal B_+&\longrightarrow
\Lambda^{1,0}T^*\St\otimes_\C E,
&
c(e_1)c(w)P_+&\longmapsto
\tfrac12(\theta^1+i\theta^2)\otimes w,\\
\Phi_-:\ p^*\mathcal B_-&\longrightarrow
\Lambda^{0,1}T^*\St\otimes_\C E,
&
c(e_1)c(w)P_-&\longmapsto
\tfrac12(\theta^1-i\theta^2)\otimes w.
\end{aligned}
\]
Indeed, under a rotation of the oriented tangent
frame, the changes in the coframe and in the
normal coefficient $w$ cancel in these expressions.
Hence, we define multiplication by $i$ on $p^*\mathcal B_\pm$
by $i\,b:=\Phi_\pm^{-1}(i\,\Phi_\pm(b))$,
so that $\Phi_\pm$ are $\C$-linear. The two components of the normal connection are then
given locally by
\[
\begin{aligned}
\partial_Ew
&=\tfrac12(\theta^1+i\theta^2)\otimes
  (\nabla_1^\perp w-J_N\nabla_2^\perp w),\\
\bar\partial_Ew
&=\tfrac12(\theta^1-i\theta^2)\otimes
  (\nabla_1^\perp w+J_N\nabla_2^\perp w).
\end{aligned}
\]
Consequently, \eqref{eq:Clocal} gives
$\Phi_+(C_+w)=\partial_Ew$ and
$\Phi_-(C_-w)=\bar\partial_Ew$,
where here $C_\pm$ denote the lifted operators.

The Hermitian dual identifies $\partial_E$ anti-linearly with
$\bar\partial_{E^*}$, so Riemann--Roch \eqref{eq:rr}, or equivalently
its Dolbeault index form \cite[Theorem 28]{Bob11}, gives the complex indices
$1-g-d_E$ and $1-g+d_E$.

The deck involution acts anti-linearly, has square $+1$, and intertwines
the lifted operators and their adjoints. To see its action on the target,
simultaneously reflect the frames by
\[
(e_1,e_2,\nu_1,\nu_2)\longmapsto(e_1,-e_2,\nu_1,-\nu_2).
\]
Both in the normal frame and in the target frame
$c(e_1)c(\nu_1)P_\pm,c(e_1)c(\nu_2)P_\pm$, the first vector is unchanged
and the second changes sign. This is complex conjugation in the
corresponding coordinates. Lemma~\ref{lem:real} therefore identifies the
real dimensions of the downstairs kernels and cokernels with the complex
dimensions upstairs. Since $1-g=\chi(\Sigma)$, we obtain
\begin{equation}\label{eq:Cindex}
\ind_\R C_++\ind_\R C_-=2(1-g)=2\chi(\Sigma).
\end{equation}

If $\Sigma$ is orientable, choose an orientation and apply the same
calculation directly to the normal complex line bundle $E:=(N,J_N)$.
Writing $h$ for its genus and $d_E$ for the degree of $E$, the complex
indices are $1-h-d_E$ and $1-h+d_E$. Passing to real indices doubles each,
and hence
\[
\ind_\R C_++\ind_\R C_-=4(1-h)=2\chi(\Sigma).
\]
Thus the sum is $2\chi(\Sigma)$ in both cases.

By Hodge theory and the Betti-number identities in
Section~\ref{subsec:forms},
\[
\dim_\R\HH=b_1(\Sigma)=
\begin{cases}
1-\chi(\Sigma),&\Sigma\text{ nonorientable},\\
2-\chi(\Sigma),&\Sigma\text{ orientable}.
\end{cases}
\]
Adding $\HH$ to the domain of $C_\pm$ increases its real Fredholm index
by $b_1(\Sigma)$. Adding the finite-rank map
$(X,w)\mapsto\bigl(c(X)-2S_X\bigr)P_\pm$ then preserves Fredholmness and the value of the index. Consequently,
\begin{equation}\label{eq:sumindex}
\ind_\R\mathcal A_++\ind_\R\mathcal A_-
=2\chi(\Sigma)+2b_1(\Sigma)
=\begin{cases}
2,&\Sigma\text{ nonorientable},\\
4,&\Sigma\text{ orientable}.
\end{cases}
\end{equation}
At least one augmented operator has positive index and therefore a
nonzero kernel. Fix such a sign and choose $(X,w)\neq(0,0)$ in that kernel.
Since $X$ is smooth and
$C_\pm w=\bigl(2S_X-c(X)\bigr)P_\pm$, elliptic regularity gives
$w\in\Gamma(N)$; see \cite[Chapter III, Section 5]{LM89}.

In summary, we obtained a nonzero smooth pair $(X,w)$ with $X^\flat$
harmonic and satisfying \eqref{eq:normalcorrection}.

\medskip\noindent\emph{Step 4: $w$ is parallel.}
For this pair $(X, w)$, \eqref{eq:Dsquare} and \eqref{eq:normpair} give
\[
0\leq\sum_{a=1}^8
\bigl(\|\mathscr D\Psi_a\|_2^2-\|\Psi_a\|_2^2\bigr)
=-4\|X\|_2^2.
\]
Thus $X=0$, while $w\neq0$ and one of $C_\pm w=0$. Each summand is nonnegative and their sum is zero. Since $\mathscr D^2-I$ is nonnegative and self-adjoint, equality implies
\[
(\mathscr D^2-I)\bigl(c(w)P_\pm\phi_a\bigr)=0
\qquad(a=1,\ldots,8).
\]
On the other hand, \eqref{eq:Dpair}, first with the chosen sign and then
with the opposite sign, gives
\begin{align*}
\mathscr D\bigl(c(w)P_\pm\phi_a\bigr) = c(w)P_\mp\phi_a, & &
\mathscr D\bigl(c(w)P_\mp\phi_a\bigr) = (C_\mp w)\phi_a+c(w)P_\pm\phi_a.
\end{align*}
It follows that $(C_\mp w)\phi_a=0$ for every $a$.
The constant sections span each real fiber, so the other $C_\mp w=0$ as well.
Both equations in \eqref{eq:Clocal} therefore have zero right-hand side.
Since $v\mapsto c(e_1)c(v)P_\pm$ is injective on $N_p$, they imply
\[
\nabla_1^\perp w-J_N\nabla_2^\perp w=0,
\qquad
\nabla_1^\perp w+J_N\nabla_2^\perp w=0.
\]
Adding and subtracting gives $\nabla^\perp w=0$.

Finally, $\Lt\leq0$ and \eqref{eq:forms} give
\[
0\leq\Qt(w)
=\int_\Sigma\bigl(|\nabla^\perp w|^2-|A_w|^2\bigr)\,d\mu_\Sigma
=-\int_\Sigma|A_w|^2\,d\mu_\Sigma.
\]
Thus $A_w=0$. The Gauss--Weingarten formula \eqref{eq:gaussweingarten}
now implies
\[
D_X^{\R^5}w=-dF(A_wX)+\nabla_X^\perp w=0
\qquad(X\in\Gamma(T\Sigma)).
\]
Connectedness of $\Sigma$ makes $w$ a fixed nonzero vector in $\R^5$.
Its normality gives $\langle F,w\rangle=0$, and hence
$F(\Sigma)\subset w^\perp\cap S^4\cong S^3$.
\end{proof}

\subsection{Superminimality in the presence of a nonzero
Dirac-harmonic section}
\begin{proposition}\label{prop:equality}
Let $F:\Sigma\to S^4$ be a closed, connected minimal immersion with $\Lt\leq0$. If $\ker \mathscr D\neq0$, then $F$ is superminimal.
\end{proposition}

\begin{proof}
Choose a nonzero smooth $\Psi\in\ker\mathscr D$ of pure chirality, i.e.
$\Omega\Psi=\varepsilon\Psi$ for some $\varepsilon\in\{+1,-1\}$.
Such a choice is possible because
$\mathscr DP_\pm=P_\mp\mathscr D$: both chiral projections preserve
$\ker\mathscr D$, and at least one projection of a nonzero section is
nonzero since $P_++P_-=I$.
The quadratic-form identity~\eqref{eq:trace} gives
\[
\sum_{a=1}^8\Qt(W_a[\Psi])=0.
\]
Since $\Lt\leq0$, each summand is nonnegative and hence zero.
Corollary~\ref{cor:zeroshiftedform} therefore implies
$\Lt W_a[\Psi]=0$ for every $a$. We will use these normal equations
to rule out a noncircular point.

\medskip\noindent\emph{Step 1: A first-order equation for $\Psi$.}
We first claim that, for every $p\in\Sigma$ and $\nu\in N_p$,
\begin{equation}\label{eq:Bzero}
c\bigl((B-sI)\nu\bigr)\Psi
-2\sum_{i=1}^2S_i\, c(\nu)\nabla_i\Psi=0.
\end{equation}
Here $B$ is defined in \eqref{eq:jacobi}, $s$ in \eqref{eq:uv}, and
$S_i=S_{e_i}$ in \eqref{eq:connection}, where $\{e_i\}$ is a local orthonormal frame of $T\Sigma$. This identity will relate the
anisotropy of $B$ to the derivatives of $\Psi$ in the next step.

To prove it, choose tangent and normal frames parallel at $p$ for
$\nabla^\Sigma$ and $\nabla^\perp$, respectively. Differentiating
\eqref{eq:fields}, using Clifford compatibility and
$\nabla_i\phi_a=S_i\phi_a$, gives
\[
\langle\nabla_i^\perp W_a[\Psi],\nu\rangle
=\langle c(\nu)\nabla_i\Psi-S_i c(\nu)\Psi,\phi_a\rangle.
\]
Let $\Delta:=-\nabla^*\nabla$ denote the nonpositive connection
Laplacian on $\SSb$. Differentiate once more, trace over $i$, and use
$(\Delta^\perp+B)W_a[\Psi]=\Lt W_a[\Psi]=0$.
Since the $\phi_a$ span each real fiber of $\SSb$, the resulting eight
scalar equations are equivalent to
\begin{equation}\label{eq:secondderivative}
0 = c(\nu)\Delta\Psi-2\sum_iS_i c(\nu)\nabla_i\Psi
  + \sum_i(S_i^2-\nabla_iS_i)c(\nu)\Psi+c(B\nu)\Psi \,.
\end{equation}
We now simplify the even-order terms. At $p$, Clifford compatibility,
Codazzi, and minimality give
\begin{align*}
\sum_i\nabla_iS_i 
= -\frac12\sum_{i,j}c(e_j)c\bigl((\nabla_{e_i}A)(e_i,e_j)\bigr) 
= -\frac12\sum_jc(e_j)c\left(\nabla_{e_j}^\perp \sum_iA(e_i,e_i)\right) = 0.
\end{align*}
Here $\nabla_iS_i$ uses the induced connection on
$\operatorname{End}_\C(\SSb)$, and $\nabla A$ uses the tangent and
normal connections. The tangent frame is parallel at $p$, so the
correction $S_{\nabla_{e_i}^\Sigma e_i}$ in the divergence vanishes there.
The computation in the proof of Proposition~\ref{prop:trace}, together
with $S_i^*=-S_i$, also gives
\[
\sum_iS_i^2=-\frac{1+s}{2}I-\frac\kappa2\Omega.
\]
Meanwhile, $\mathscr D\Psi=0$ and the Weitzenb\"ock identity \eqref{eq:weitzenbock} imply
\[
\Delta\Psi=\left(\frac K2I-\frac\kappa2\Omega\right)\Psi.
\]
Using $K=1-s$ and $c(\nu)\Omega=-\Omega c(\nu)$, substitution these equations into \eqref{eq:secondderivative} proves \eqref{eq:Bzero}.

\medskip\noindent\emph{Step 2.}
At every point where the curvature ellipse is noncircular, we claim that
\begin{equation}\label{eq:gradientnorm}
|\nabla\Psi|^2=\frac{s+\varepsilon\kappa}{2}|\Psi|^2,
\end{equation}
where $s$ and $\kappa$ are defined in \eqref{eq:uv}.

To prove this, choose orthonormal tangent and normal frames along the
ellipse axes at such a point, compatible with the prescribed combined
orientation, and let
\[
\alpha:=\langle A_{11},\nu_1\rangle,
\qquad\beta:=\langle A_{12},\nu_2\rangle.
\]
Then $u=A_{11}=\alpha\nu_1$, $v=A_{12}=\beta\nu_2$, and
$\alpha^2\neq\beta^2$. In particular,
\[
s=\alpha^2+\beta^2,\qquad\kappa=2\alpha\beta,
\qquad(B-sI)\nu_1=(\alpha^2-\beta^2)\nu_1.
\]
The equation $\mathscr D\Psi=0$ gives
$\nabla_2\Psi=-c(e_1)c(e_2)\nabla_1\Psi$.
Moreover, \eqref{eq:connection} and minimality give
\begin{align*}
2S_1&=c(e_1)-c(e_1)c(u)-c(e_2)c(v),\\
2S_2&=c(e_2)-c(e_1)c(v)+c(e_2)c(u),
\end{align*}
so $S_2c(e_1)c(e_2)=c(e_1)-S_1$.
Since $c(\nu)$ commutes with $c(e_1)c(e_2)$, it follows that
\[
\sum_iS_i c(\nu)\nabla_i\Psi
=\bigl(S_1-S_2c(e_1)c(e_2)\bigr)c(\nu)\nabla_1\Psi
=\bigl(2S_1-c(e_1)\bigr)c(\nu)\nabla_1\Psi.
\]
Substituting into \eqref{eq:Bzero} with $\nu=\nu_1$ and multiplying
on the left by $c(\nu_1)$ gives
\[
2\bigl(\alpha c(e_1)c(\nu_1)-\beta c(e_2)c(\nu_2)\bigr)\nabla_1\Psi
=(\alpha^2-\beta^2)\Psi.
\]
Taking norms on both sides yields
\begin{equation}\label{eq:derivativenormfactor}
4(s-\varepsilon\kappa)|\nabla_1\Psi|^2
 = (\alpha^2-\beta^2)^2|\Psi|^2
 = (s-\varepsilon\kappa)(s+\varepsilon\kappa)|\Psi|^2.
\end{equation}
Here we use the Clifford identity
\begin{align*}
 \bigl(\alpha c(e_1)c(\nu_1)-\beta c(e_2)c(\nu_2)\bigr)^*
 \bigl(\alpha c(e_1)c(\nu_1)-\beta c(e_2)c(\nu_2)\bigr)
 = sI-\kappa\Omega \,;
\end{align*}
and the fact that $\Omega$ is parallel, so $\nabla_1\Psi$ has the same
chirality $\varepsilon$ as $\Psi$.
At a noncircular point $s>|\kappa|$, by \eqref{eq:circularitycriteria}.
Thus $s-\varepsilon\kappa>0$. Canceling this factor in
\eqref{eq:derivativenormfactor} and using
$|\nabla\Psi|^2=2|\nabla_1\Psi|^2$ proves \eqref{eq:gradientnorm}.

\medskip\noindent\emph{Step 3: circularity everywhere.}
Suppose there is a noncircular point. Apply Lemma~\ref{lem:quartic}
on $\Sigma$ if it is orientable, and on its connected orientation cover
otherwise. The holomorphic quartic is nonzero at that point or either
of its lifts, so it is not identically zero. Its zeros, which are
precisely the circular points, are therefore isolated. The noncircular
points are dense in $\Sigma$; in the nonorientable case this follows
downstairs as well because the covering map is a local diffeomorphism.
Consequently \eqref{eq:gradientnorm} extends everywhere by continuity.

We can now integrate on the closed surface. Since $\mathscr D\Psi=0$
and $\Omega\Psi=\varepsilon\Psi$, equations \eqref{eq:weitzenbock}
and \eqref{eq:gradientnorm} give
\begin{align*}
0= \|\mathscr D\Psi\|_2^2
 = \int_\Sigma\left(|\nabla\Psi|^2 +\frac{K-\varepsilon\kappa}{2}|\Psi|^2\right)\,d\mu_\Sigma 
 = \frac12\int_\Sigma(s+K)|\Psi|^2\,d\mu_\Sigma
=\frac12\|\Psi\|_2^2.
\end{align*}
The last equality uses $K+s=1$ from \eqref{eq:gaussricci} and contradicts
the choice $\Psi\neq0$. Thus the curvature ellipse is circular
everywhere, and so $F$ is superminimal.
\end{proof}

\begin{proof}[Proof of Theorem~\ref{thm:circular}]
Thanks to Propositions~\ref{prop:paired} and~\ref{prop:equality} we are only left with justifying the last assertion. Suppose first that $\Sigma$ is nonorientable and $F$ is an embedding. A closed nonorientable surface cannot embed in $\R^3$; see \cite[Corollary 3.46 and the following paragraph, p.~256]{Hat02}. Removing a point outside the image gives the same obstruction in $S^3$. Thus $F(\Sigma)$ cannot lie in a great $S^3$. For a nonorientable embedding, Proposition~\ref{prop:paired} must therefore yield a nonzero Dirac-harmonic section, and Proposition~\ref{prop:equality} gives superminimality.

If $\Sigma$ is orientable and $F(\Sigma)$ lies in a
great $S^3$, the immersion admits a global unit normal
$\nu$ within that three-sphere. Exploiting the assumption $\Lt\leq0$, the scalar second
variation formula with constant (unit) speed then gives 
$A\equiv0$, so $F$ is totally geodesic and thus
superminimal. Thereby the proof is complete.
\end{proof}

The next statement follows as a special case of the recent result by Bernstein and Ketover \cite{BK26a} cited in the introduction. However, in ambient dimension four we can easily embed the justification of the improved Morse index bound in our discussion.

\begin{corollary}[The orientable case]\label{cor:orientableindex}
A closed, connected, orientable minimal immersion into $S^4$ that is not totally geodesic has index at least six.
\end{corollary}
\begin{proof}
Proposition~\ref{prop:five} gives index at least five. If the index were five, it would also give $\Lt\leq0$ and $\dim_\R\ker\Lt=5$. The preceding proof and Theorem~\ref{thm:circular} imply superminimality; in fact, the great-three-sphere alternative would force the immersion to be totally geodesic.
Choose an orientation of $\Sigma$ and the compatible normal orientation. Then $J_N:N\to N$ is a global parallel orthogonal complex structure. By~\eqref{eq:circularitycriteria}, superminimality gives $B=sI$, so both $\Delta^\perp$ and $B$ commute with $J_N$. Consequently $\Lt J_N=J_N\Lt$, and $\ker\Lt$ is invariant under $J_N$. This is the standard complex-eigenspace observation for superminimal surfaces; see \cite[Section 3]{MU97}. Thus $\ker\Lt$ is invariant under the complex
structure $J_N$, so it has even real dimension,
contradicting $\dim_\R\ker\Lt=5$.
\end{proof}

\section{The superminimal case and completion of the proof}\label{sec:degree}
We use two classical results for superminimal surfaces, stated here in the normalization of the unit round $S^4$.

\medskip\noindent\textbf{Holomorphic normal kernel (Montiel--Urbano).}
Let $f:M\to S^4$ be a closed, connected, oriented superminimal immersion of genus $g$ that is not totally geodesic. For a connected oriented superminimal surface, either $\kappa=s$ everywhere or $\kappa=-s$
everywhere; see the argument immediately preceding
Theorem~6.5 in \cite[p.~611]{ES85}.
Reversing the ambient orientation changes the sign
of $\kappa$ and leaves $s$ unchanged.
We henceforth choose that orientation so that
$\kappa=s$ on $M$.
Let $N\to M$ be the normal bundle of $f$, and let $\Lt$
be its shifted normal Jacobi operator, as defined in
\eqref{eq:shifted}. Equip the complex line bundle
$E=(N,J_N)\to M$ with the holomorphic structure defined
by the $(0,1)$ part of the normal connection. Then
\begin{equation}\label{eq:classicalkernel}
\ker\Lt=H^0(E),\qquad
h^0(E)=\frac{\Area(M)}{2\pi}+g-1.
\end{equation}
Here $H^0(E)$ denotes the space of global holomorphic sections
of $E$, identified with real normal fields through its
underlying real bundle, and $h^0(E):=\dim_\C H^0(E)$.

These are the holomorphic kernel identification and
Riemann--Roch calculation of Montiel--Urbano,
see Theorem 3 in \cite{MU97} and its proof; see also
\cite[Lemma 4.1 and Theorem 4.2]{KW24}.
Since $-L=-\Lt-2I$, this kernel is precisely the eigenspace
of the second-variation operator $-L$ with eigenvalue $-2$.

\medskip\noindent\textbf{Twistor correspondence and degree bound.}
Such an immersion has a horizontal holomorphic lift $\widehat f:M\to\CP^3$ for the twistor projection $\pi:\CP^3\to S^4$, after the same choice of ambient orientation. Its image is full, meaning that it is not contained in a complex projective hyperplane of $\CP^3$, and
\begin{equation}\label{eq:areadegree}
\Area(M)=4\pi d,\qquad d:=\deg\widehat f^*\mathcal O(1)\geq3.
\end{equation}
Here $\mathcal O(1)$ is the hyperplane line bundle on $\CP^3$, and horizontal means that the differential of $\widehat f$ takes values in the horizontal distribution of the twistor fibration. We use the correspondence of Bryant \cite{Bry82} and Eells--Salamon \cite[Theorem 6.5]{ES85}; fullness and the area--degree relation are recalled in \cite[Proposition 9.4 and the subsequent paragraph]{ES85}. The resulting bound $d\geq3$ is also recorded in \cite[Section 4]{MU97} and \cite[proof of Theorem 4.3]{KW24}.

In the application below, superminimality follows from
Theorem~\ref{thm:circular} and is preserved upon lifting
the immersion to its orientation double cover.

\begin{proposition}\label{prop:degree}
Let $F:\Sigma\to S^4$ be a closed, connected, nonorientable superminimal immersion. Suppose $\dim_\R\ker\Lt=5$. Then the orientation cover has genus zero, its horizontal holomorphic twistor lift has degree three, and $\Sigma$ is diffeomorphic to $\RP^2$.
\end{proposition}
\begin{proof}
Let $p:\St\to\Sigma$ be the connected orientation double cover, of genus $g$, and let $\widetilde F=F\circ p$. The immersion cannot be totally geodesic, since a closed connected totally geodesic surface covers a great $S^2$ and is therefore orientable, as recalled in Section~\ref{subsec:simons}. The holomorphic normal kernel identification
\eqref{eq:classicalkernel} and the twistor lift with
the area--degree relation \eqref{eq:areadegree}
therefore apply to $\widetilde F$, with the ambient
orientation chosen so that $\kappa=s$ as explained
at the beginning of this section.

Let $E=(p^*N,J_N)\to\St$ be its holomorphic normal line, with the sign chosen as above. 
Let $\tau:\St\to\St$ be the nontrivial deck
transformation of the orientation cover.
It reverses the orientation of $\St$.
Since $\widetilde F\circ\tau=\widetilde F$, the normal
fibers at $x$ and $\tau(x)$ are naturally identified
as subspaces of the same ambient tangent space.
This identification reverses their induced normal
orientations, because the ambient orientation is
fixed whereas the tangent orientation is reversed.

Under the canonical identification
$(p^*N)_x=N_{p(x)}=(p^*N)_{\tau(x)}$, reversal of the
normal orientation gives $J_N(\tau(x))=-J_N(x)$.
Consequently, the map
\[
T:\Gamma(p^*N)\longrightarrow\Gamma(p^*N),
\qquad (TV)(x):=V(\tau(x)),
\]
where the value on the right is viewed in $(p^*N)_x$
through this identification, is an anti-linear
involution. A section $V\in\Gamma(p^*N)$ satisfies $TV=V$
if and only if there exists a section
$W\in\Gamma(N)$ such that
\[
V(x)=W(p(x))
\qquad(x\in\St),
\]
where $N\to\Sigma$ is the normal bundle of the
original immersion $F$.

Denote by
$\widehat{\Lt}:\Gamma(p^*N)\longrightarrow\Gamma(p^*N)$
the lift of the shifted normal Jacobi operator
$\Lt:\Gamma(N)\to\Gamma(N)$.
Because its coefficients and connection are pulled
back from $\Sigma$, it satisfies
\[
\widehat{\Lt}(p^*W)=p^*(\Lt W)
\qquad(W\in\Gamma(N)),
\qquad
\widehat{\Lt}\,T=T\,\widehat{\Lt}.
\]
In particular, $T$ preserves $\ker\widehat{\Lt}$.
By \eqref{eq:classicalkernel}, this kernel identifies
with $H^0(E)$, so $T$ restricts to an anti-linear
involution on $H^0(E)$. Define its fixed subspace by
$H^0(E)^T:=\{V\in H^0(E):TV=V\}$.
The pullback map restricts to an $\R$-linear
isomorphism
\[
p^*:\ker\Lt\longrightarrow H^0(E)^T.
\]
Indeed, an invariant section has the form $p^*W$,
and the preceding operator identity shows that
$\widehat{\Lt}(p^*W)=0$ if and only if $\Lt W=0$.

The dimension identity for an anti-linear involution,
proved in Lemma~\ref{lem:real}, now gives
\[
\dim_\R H^0(E)^T=\dim_\C H^0(E)=h^0(E).
\]
Using the hypothesis $\dim_\R\ker\Lt=5$ and
\eqref{eq:classicalkernel}, we conclude that
\begin{equation}\label{eq:areakernel}
5=\dim_\R\ker\Lt
=\dim_\R H^0(E)^T
=h^0(E)
=\frac{\Area(\St)}{2\pi}+g-1.
\end{equation}
Let $d$ be the degree of the horizontal holomorphic
twistor lift of $\widetilde F$.
Equation~\eqref{eq:areadegree} gives
$\Area(\St)=4\pi d$ and $d\geq3$. Hence $2d+g=6$ with $d\geq3$ and obviously $g\geq0,$
which forces $d=3$ and $g=0$.
Since $\chi(\Sigma)=\tfrac12\chi(\St)=1,$
the classification of closed surfaces implies that
$\Sigma$ is diffeomorphic to $\RP^2$.
\end{proof}

To finish we use Karpukhin's index characterization \cite{Kar21}, in the form stated in \cite[Theorem 4.8]{KW24}: a minimal immersion of $\RP^2$ into $S^n$, $n\geq4$, has index at least $n+1$, and equality holds precisely for the Veronese plane in a great $S^4$. Thus once the topology has been determined, no separate classification of horizontal cubics is needed.

\begin{proof}[Proof of Theorem~\ref{thm:main}]
Proposition~\ref{prop:five} gives $\Ind(\Sigma)\geq5$, since the immersion is not totally geodesic. Corollary~\ref{cor:orientableindex} gives the stronger bound in the orientable case. If $\Ind(\Sigma)=5$, then $\Sigma$ is nonorientable, $\Lt\leq0$, and $\dim_\R\ker\Lt=5$. Theorem~\ref{thm:circular} gives superminimality, and Proposition~\ref{prop:degree} gives $\Sigma\cong\RP^2$. Karpukhin's characterization with $n=4$ now implies that $F(\Sigma)=R\Pi$ for some $R\in O(5)$.
\end{proof}

\vspace{1mm}
 
\section*{Acknowledgements}
This project has received funding from the European Research Council (ERC) under the European Union’s Horizon 2020 research and innovation programme (grant agreement No. 947923). The second author was partially supported by Simons Foundation Travel Support for Mathematicians 2026-2031.

\vspace{1mm}
\stoptoc
\section*{AI Declaration Statement}
This project naturally stemmed from the previous works of the authors about minimal submanifolds with isolated singularities \cite{CLW24, CLW26}. 
Over the past several months, the authors obtained a series of partial results toward a complete classification of minimal surfaces of Morse index five, including a treatment of the superminimal case and the exclusion of orientable surfaces as well as of immersions whose normal bundles have large Euler number. 
Hence, after extensive interactions relying upon such inputs, a LLM (Chat GPT-6 Astra) suggested the approach that eventually led to Theorem \ref{thm:circular}. The authors independently verified the correctness of the strategy and turned that into a complete proof. This article contains no AI-generated text. 

\resumetoc

\printbibliography

\end{document}